\documentclass[11pt]{amsart}

\usepackage{amsfonts}
\usepackage{amsmath}
\usepackage{amssymb, esint}
\usepackage{amsthm,color}
\usepackage{enumerate}
\usepackage{mathtools}
\usepackage{cite}

\newcommand{\R}{{\mathbb R}}

\newcommand{\N}{{\mathbb N}}

\newcommand{\cS}{{\mathcal S}}

\newcommand{\cP}{{\mathcal P}}

\def\0{{\mathbf 0}}

\newcommand{\e}{\varepsilon}
\newcommand{\vp}{\varphi}
\newcommand{\osc}{\operatornamewithlimits{osc}}

\newcommand{\supp}{\operatorname{spt}}

\newcommand{\dist}{\operatorname{dist}}

\newcommand{\tr}{\operatorname{tr}}

\newcommand\norm[1]{\Arrowvert {#1}\Arrowvert}
\newcommand\inn[2]{\langle {#1,#2} \rangle}

\def\mean#1{\mathchoice%
          {\mathop{\kern 0.2em\vrule width 0.6em height 0.69678ex depth -0.58065ex
                  \kern -0.8em \intop}\nolimits_{\kern -0.4em#1}}%
          {\mathop{\kern 0.1em\vrule width 0.5em height 0.69678ex depth -0.60387ex
                  \kern -0.6em \intop}\nolimits_{#1}}%
          {\mathop{\kern 0.1em\vrule width 0.5em height 0.69678ex
              depth -0.60387ex
                  \kern -0.6em \intop}\nolimits_{#1}}%
          {\mathop{\kern 0.1em\vrule width 0.5em height 0.69678ex depth -0.60387ex
                  \kern -0.6em \intop}\nolimits_{#1}}}

\theoremstyle{plain}
\newtheorem{thm}{Theorem}[section]

\newtheorem{cor}[thm]{Corollary}
\newtheorem{lem}[thm]{Lemma}
\newtheorem{prop}[thm]{Proposition}

\newtheorem{defn}[thm]{Definition}

\theoremstyle{rem}

\newtheorem{case}{Case}

\newtheorem*{claim*}{Claim}

\newtheorem{rem}[thm]{Remark}

\numberwithin{equation}{section}

\title[]{A transmission problem with $(p,q)$-Laplacian}

\author{Maria Colombo}
\address{Institute of Mathematics, EPFL SB, Station 8, CH-1015 Lausanne, Switzerland
maria.colombo@epfl.ch
}
\email{maria.colombo@epfl.ch}

\author{Sunghan Kim}
\address{Department of Mathematics, Royal Institute of Technology, 100 44 Stockholm, Sweden}
\email{sunghan@kth.se}

\author{Henrik Shahgholian}
\address{Department of Mathematics, Royal Institute of Technology, 100 44 Stockholm, Sweden}
\email{henriksh@kth.se}

\thanks{M. Colombo thanks the SNSF Grant 182565. S. Kim was supported by postdoctoral fellowship from Knut and Alice Wallenberg Foundation. H. Shahgholian was supported in part by Swedish Research Council.   The authors would like to thank John Andersson  for sharing ideas and having fruitful discussions on Section \ref{section:holder}.}

\begin{document}

\maketitle

\begin{abstract}
In this paper we consider  the so-called double-phase problem where the phase transition takes place across
 the zero level "surface" of the minimizer of the functional 
$$
J(v,\Omega) = \int_\Omega \left( |D v^+|^p +  |D v^-|^q \right) dx.
$$
We prove that the minimizer exists, and is  H\"older regular. From here, using  an intrinsic variation,
 one can prove  a weak formulation of the free boundary condition  across the zero level surface, that formally 
can be represented as
 $$
 (q-1)|D u^-|^q = (p-1) |D u^+|^p  \quad \hbox{on } \partial \{u > 0\}.
$$
 We prove that the free boundary is $C^{1,\alpha}$ a.e. 
 with respect to the measure $\Delta_p u^+$, whose support is of $\sigma$-finite $(n-1)$-dimensional Hausdorff measure.
\end{abstract}

\tableofcontents


\section{Introduction}\label{section:intro}

\subsection{Background}
Since the  works of Zhikov  \cite{Zhikov1},  \cite{Zhikov2},  there has been a surge in studying variational problems with  non-standard growth of the type
\begin{equation}\label{eq:main}
\int a(x,u) |D u|^{\Lambda(x,u)} + b(x,u)
\end{equation}
with many variations and conditions on $\Lambda(x,u)$. Although Zhikov's motivation concerned homogenization and the so-called Lavrentiev's phenomenon, the model problem seems to have other applications in electrodynamics of moving media in material sciences such as  electrorheological fluids;  see \cite{Ruzicka} and the references therein. For the development of the mathematical theory we refer to \cite{Colombo-Mingione-1, Colombo-Mingione-2, BCM18}, \cite{Acerbi-Fusco}, \cite{Chipot}, \cite{LW17,LW19}, just to have mentioned a few. 

Our prime goal in this paper is to set this topic into a new perspective in terms of material dependent conductivity that 
may give discontinuity in the conduction across a level set (which a priori is unkonw).
 Such problems for the case of standard-growth (Dirichlet energy) were studied in \cite{Kiam-Kim-Sh1, Kiam-Kim-Sh2}, where complete results were obtained.
We shall thus consider the case of $\Lambda(x,u) = p\chi_{\{u > 0\}} + q \chi_{\{ u < 0\}}$ in \eqref{eq:main}, but our approach maybe applied to more general setting such as $p,q$ may  also depend on $x$, and $u$. 

To formalize the setting we let $n\geq 2$ be the space dimension, $p,q\in(1,\infty)$ with $p\neq q$, and $\Omega\subset\R^n$ be a domain. Our model problem is given by the functional 
\begin{equation}\label{eq:J}
J(v,\Omega) \equiv J_{p,q}(v,\Omega) := \int_\Omega ( |Dv^+|^p + |Dv^-|^q ) \, dx,
\end{equation}
where $v^+ := \max\{v,0\}$, $v^- := -\min\{v,0\}$. We should  remark that the functional $J$ is not convex, and hence uniqueness may fail in general. 

\begin{defn}\label{defn:min}
We shall call $u\in W_{loc}^{1,p\wedge q}(\Omega)$ a local minimizer, of $J$ defined in \eqref{eq:J}, if $J(u,\Omega')<\infty$ for any bounded domain $\Omega'\Subset \Omega$, and $J(u,\supp \phi) \leq J( u + \phi , \supp \phi)$ for any $\phi\in W_{loc}^{1,p\wedge q}(\Omega)$ satisfying $\supp\phi\subset\Omega$. 
\end{defn}

The definition above implies that if $u \in W_{loc}^{1,p\wedge q}(\Omega)$ is a local minimizer, then $u^+ \in W_{loc}^{1,p}(\Omega)$ and $u^- \in W_{loc}^{1,q}(\Omega)$. Since we shall analyse local properties of local minimizers, we shall only treat local minimizers $u\in W^{1,p\wedge q}(\Omega)$ satisfying $J(u,\Omega) < \infty$. We shall simply call $u\in W_{loc}^{1,p\wedge q}(\Omega)$ a local minimizer, whenever the functional $J$ is well-understood in the context. 

We shall also say that $u\in W_{loc}^{1,p}(\Omega)$ is a weak $p$-subsolution (resp. $p$-supersolution), if
$$
\int_\Omega |Du|^{p-2} Du \cdot D \phi\,dx \leq \text{(resp. $\geq$) } 0, 
$$
for any $\phi\in C_0^\infty(\Omega)$. We shall call $u\in W_{loc}^{1,p}(\Omega)$ a $p$-harmonic function, if $u$ is simultaneously a weak $p$-subsolution and a weak $p$-supersolution. For basic properties of $p$-harmonic functions, see \cite{L}. 

\subsection{Main results}

Before closing this introduction we state our main results. The existence theory, Proposition \ref{prop:exist}, follows in a standard way, and we have included only for the benefit of readers. The first result is about the H\"older regularity of local minimizers of the functional $J$. 
\begin{thm}\label{thm:holder}
Let $u \in W^{1,p\wedge q}(\Omega)$ be a local minimizer of the functional $J$. Then $u\in C_{loc}^{0,\alpha}(\Omega)$ for some $\alpha\in(0,1)$ depending only on $n$, $p$ and $q$. Moreover, for every compact set $K\subset\Omega$, 
\begin{equation}\label{eq:Ca}
[u^\pm]_{C^{0,\alpha}(K)} \leq \frac{C}{(\dist(K,\partial\Omega))^{\alpha+\frac{n}{p_\pm}}}\norm{u^\pm}_{L^{p_\pm}(\Omega)},
\end{equation}
where $p_+ = p$, $p_- = q$ and $C>0$ depends only on $n$, $p$ and $q$. 
\end{thm}

The main difficulties arise from the fact that each phase ($u^+$, $u^-$) scales differently (see Lemma \ref{lem:scale}), as the energy functional exhibits jump discontinuity in the power, across the zero-level set. The key is to obtain a geometric decay of each phase, regardless of the behavior of the other. The proof is based on several technical results, in Section \ref{section:basic}, that may be of independent interest in general. 

\begin{rem}
Our proof is based on the universal decay of each phase around the zero-level set (Proposition \ref{prop:holder}). As a contraposition (due to the De Giorgi lemma), we prove a weak-Harnack-type estimate (Proposition \ref{prop:holder-re}) that if the negative (positive) phase is small in terms of measure, then the minimizer is strictly positive (negative) in the interior. We would like to emphasize that the classical approach fails in our setting, due to the presence of free boundaries. Especially, our minimizers admit subsolution properties, but lack supersolution properties (in the sense of the De Giorgi class \cite{Giu}), which makes our analysis difficult and non-standard. 
\end{rem}

The second result is on the $C^{1,\alpha}$-regularity of the zero-level surface (also referred to as the free boundary) of local minimizers, provided that it is universally flat. 

\begin{thm}\label{thm:flat}
Let $u \in W^{1,p\wedge q}(\Omega)$ be a local minimizer of the functional $J$, $z\in\partial\{u>0\}\cap\Omega$ be given, and $0<r<\dist(z,\partial\Omega)$. Then there exist some $\sigma\in(0,1)$, $\e>0$ and $\rho>0$, depending at most on $n$, $p$ and $q$, such that if $\partial\{u>0\}\cap B_r(z) \subset \{x\in\R^n: |\inn{x-z}{\nu}| \leq r\e\}$ for some direction $\nu\in\partial B_1$, then $\partial\{u>0\}\cap B_{\rho r}(z)$ is a $C^{1,\sigma}$-graph in direction $\nu$. 
\end{thm}
This is a standard result for free boundary problems of Bernoulli type. The key here is to verify that local minimizers are viscosity solutions (in the sense of Definition \ref{defn:visc}). Once this is done (Proposition \ref{prop:visc}), the rest of the argument follows essentially from the work of Lewis-Nyström \cite{LN1}, \cite{LN2}. 

Our final result concerns structure theorem for the free boundary, in terms of the measure $\Delta_p u^+$:
\begin{thm}\label{thm:nonflat}
Let $u\in W^{1,p\wedge q}(\Omega)$ be a local minimizer of the functional $J$. Then $\Delta_p u^+$ is a nonnegative Radon measure, with $\supp \Delta_p u^+ \subset \partial\{u>0\}\cap\Omega$. Moreover, for $\Delta_p u^+$-a.e.\ point $z\in\partial\{u>0\}\cap\Omega$, there exists some $r=r_z \in(0,\frac{1}{2}\dist(z,\partial\Omega))$ such that $\partial\{u>0\}\cap B_r(z)$ is a $C^{1,\sigma}$-graph, where $\sigma\in(0,1)$  depend only on $n$, $p$ and $q$.
 Furthermore, $\supp(\Delta_p u^+)$ is of $\sigma$-finite $H^{n-1}$-measure. 
\end{thm}

The last theorem is related to the work of Andersson-Mikayelyan \cite{AM}, which deals with the uniformly elliptic problems featuring jump discontinuity in the conductivity matrix across the free boundary. Our argument here exhibits some new features, as each phase scales differently, and thus is technically more involved.  

\begin{rem}
It should be stressed that $\supp(\Delta_p u^+)$ being of $\sigma$-finite $H^{n-1}$-measure is (much) weaker than $\partial\{u > 0\}$ having finite $H^{n-1}$-measure. The latter is a highly challenging problem even for $p$-harmonic functions with $p\neq 2$. 
\end{rem}

We would like to remark that the H\"older regularity theory (Theorem \ref{thm:holder}) can be extended to local minimizers of a more general functional
$$
\int a_p(x,Du^+) + a_q(x, Du^-) + f(x,u), 
$$
without major modifications, provided that $a_p$ (and $a_q$) is bounded measurable and has standard $p$-(resp. $q$-)growth, and $f$ bounded in $(x,u)$. Note that when $f(x,u) = g(x) \chi_{\{u>0\}}$, this corresponds to classical two phase Bernoulli-type problem \cite{ACF}. Thus, it is also possible to prove the $C^{1,\alpha}$-regularity of flat free boundaries (Theorem \ref{thm:flat}) for local minimizers of the above extended functional. We shall not do it here, and hope to come back to this issue in the nearest future. 
 
As a final remark, let us address that the Lipschitz regularity of local minimizers is left open. This is a very challenging issue, and it would involve some sophisticated analysis on global minimizers (i.e., the minimizers on the entire space). We would like to invite the interested reader to explore more in this direction.

\subsection{Outline} 

In Section \ref{section:exist}, we prove the existence of minimizers subject to a prescribed boundary condition. In Section \ref{section:basic}, we present some basic properties and technical lemmas that will be used extensively throughout the paper. Section \ref{section:holder} is devoted to the study of H\"older regularity of local minimizers, and we prove our first result, Theorem \ref{thm:holder}. Section \ref{section:flat} concerns the regularity of flat free boundaries, and we prove our second result, Theorem \ref{thm:flat}. In Section \ref{section:nonflat}, we study the structure of the free boundary, and prove the final result, Theorem \ref{thm:nonflat}. In Appendix \ref{appendix:flat}, we present supplementary details for some of our argument, for the reader's benefit.



\section{Existence}\label{section:exist}

This section is devoted to the existence of minimizers. The proof is standard, but we contain it for the sake of completeness.

\begin{prop}\label{prop:exist} 
Let $\Omega\subset\R^n$ be a bounded domain, and let $\vp \in W^{1,p\wedge q}(\Omega)$ be such that $J(\vp,\Omega)<\infty$. Then there exists $u\in \vp + W_0^{1,p\wedge q}(\Omega)$ such that $J(u,\Omega) < \infty$, and $J(u,\Omega) \leq J(v,\Omega)$ for any $v \in \vp + W_0^{1,p\wedge q}(\Omega)$.
\end{prop}

\begin{proof} 
Let us assume the case $p>q$ only; the other case, $q>p$, follows with the same argument, so will be omitted. Here $C$ denotes a constant independent of $k$, and it may vary at each occurrence. Since $J(\vp,\Omega)<\infty$, and $J(v,\Omega)\geq 0$ for any $v\in \vp + W_0^{1,q}(\Omega)$, it suffices to prove the weak lower semicontinuity of the functional. Let $A$ be the infimum of $J(v,\Omega)$ over all $v\in \vp + W_0^{1,q}(\Omega)$. Take a sequence $\{u_k\}_{k=1}^\infty\subset W^{1,q}(\Omega)$ such that 
\begin{equation}\label{eq:exist-1}
\lim_{k\to\infty} J(u_k,\Omega) = A. 
\end{equation} 
Due to \eqref{eq:exist-1}, we may assume without loss of generality that 
\begin{equation}\label{eq:exist-2}
J(u_k,\Omega) \leq C.
\end{equation}
for all sufficiently large $k$. 

Since $q<p$ and $\Omega$ is a bounded domain, Jensen's inequality implies that  
$$
\begin{aligned}
\int_\Omega |D u_k|^q\,dx &\leq |\Omega|^{q/p} \left(\int_\Omega |D u_k^+|^p\,dx\right)^{q/p} + \int_\Omega |D u_k^-|^q\,dx \\
&\leq C(J(u_k,\Omega) + 1)  \leq C,
\end{aligned} 
$$
where the last inequality holds for all large $k$'s, due to \eqref{eq:exist-2}. Recalling that $u_k - \vp \in W_0^{1,q}(\Omega)$, this also implies along with Poincar\'e inequality that
\begin{equation}\label{eq:exist-3}
\begin{aligned}
\int_\Omega |u_k|^q\,dx &\leq 2^q \left( \int_\Omega |\vp|^q\,dx + \int_\Omega |u_k - \vp|^q\,dx \right) \\
&\leq 2^q \int_\Omega |\vp|^q\,dx + C\int_\Omega |D u_k - \vp|^q\,dx  \\
&\leq C \int_\Omega ( |\vp|^q + |D \vp|^q)\,dx + C\int_\Omega |D u_k|^q \,dx \leq  C.
\end{aligned}
\end{equation}
This shows that $\{u_k\}_{k=1}^\infty$ is a bounded sequence in $W^{1,q}(\Omega)$. 

Thus, there exist certain $u\in W^{1,q}(\Omega)$ and a subsequence of $\{u_k\}_{k=1}^\infty$, which we shall denote with the same index for the notational convenience, such that $u_k \to u$ weakly in $W^{1,q}(\Omega)$. Extracting a further subsequence if necessary, we know that $D u_k^\pm \to D u^\pm$ weakly in $L^q(\Omega)$, and $u_k^\pm \to u^\pm$ strongly in $L^q(\Omega)$. 

It is also straightforward from \eqref{eq:exist-2} that
$$
\int_\Omega |D u_k^+|^p\,dx \leq C, 
$$
and hence, by $u_k^+ - \vp^+ \in W_0^{1,p}(\Omega)$, we can follow the lines of \eqref{eq:exist-3} and deduce that
$$
\int_\Omega |u_k^+|^p\,dx \leq C. 
$$
Thus, $\{u_k^+\}_{k=1}^\infty$ is also a bounded sequence in $W^{1,p}(\Omega)$, from which it follows that $D u_k^+ \to \hat u$ weakly in $L^p(\Omega)$ and $u_k^+ \to \hat u$ strongly in $L^p(\Omega)$, after extracting a further subsequence if necessary. However, since we have already observed $u_k^+ \to u^+$ strongly in $L^q(\Omega)$, which also implies $u_k^+ \to u^+$ a.e.\ in $\Omega$ along another subsequence, we can ensure that $\hat u = u^+$ a.e.\ on $\Omega$. This combined with $D u_k^+ \to D \hat u$ weakly in $L^p(\Omega)$ proves that $u_k^+ \to u^+$ strongly in $L^p(\Omega)$ and $D u_k^+ \to D u^+$ weakly in $L^p(\Omega)$.

Thanks to the weak convergence of $u_k^+ \to u^+$ in $W^{1,p}(\Omega)$ and $u_k^- \to u^-$ in $W^{1,q}(\Omega)$, we conclude that 
$$
J(u,\Omega) \leq \liminf_{k\to\infty} J(u_k,\Omega), 
$$
which combined with \eqref{eq:exist-2} proves the minimality of $J(u,\Omega)$ among all $v\in \vp + W_0^{1,q}(\Omega)$. Hence, the proof is complete.
\end{proof}


\section{Basic Properties}\label{section:basic}

The main purpose of this section is to present some basic properties of local minimizers $u$ of the functional $J$ as in \eqref{eq:J}; we shall simply call $u$ a local minimizer without specifying the functional, unless there arises any ambiguity. Throughout this section, $\Omega\subset\R^n$ will be a bounded domain, and $c$, $C$ will denote generic constants depending only on $n$, $p$ and $q$, unless stated otherwise; these constants may vary at each occurrence. 

We remark that many properties are symmetric in between $u^+$ and $u^-$, and for this reason we shall only present the statement and the proof for $u^+$. As a matter of fact, the argument does not distinguish the case $p>q$ nor $p<q$, so one may derive the corresponding assertion for $u^-$ by noting the following lemma. We omit the obvious proof. 

\begin{lem}\label{lem:basic}
Let $u \in W^{1,p\wedge q}(\Omega)$ be a local minimizer of the functional $J(\cdot;a_p,a_q)$. Then $(-u)$ is a local minimizer of the functional $J(\cdot;a_q,a_p)$. 
\end{lem}

It is noteworthy that the positive phase scales differently from the negative phase. 

\begin{lem}\label{lem:scale}
Let $u\in W^{1,p\wedge q}(\Omega)$ be a local minimizer, $B\subset\Omega$ be a ball and $S>0$ be a number. Let $\Omega_B\subset\R^n$ be the image of $\Omega$ via the dilation that maps $B$ to $B_1$, and let $x_B$ and $\rho$ be the centre and respectively the radius of $B$. Define $w:\Omega_B\to\R$ by 
$$
w(x) = \frac{u^+(\rho x + x_B)}{S} - \frac{u^-(\rho x + x_B)}{\rho^{1-\frac{p}{q}}S^{\frac{p}{q}}}.
$$
Then $w\in W^{1,p\wedge }(\Omega_B)$ and it is also a local minimizer.
\end{lem} 

\begin{rem}\label{rem:scale}
In most cases, $S = \norm{u^+}_{L^\infty(B)}$ (resp. $S = \norm{u^-}_{L^\infty(B)}$ etc.), so that $\norm{w^+}_{L^\infty(B_1)} = 1$ (resp. $\norm{w^-}_{L^\infty(B_1)} = 1$ etc.). The reason that we do not take, e.g., $S = \norm{u}_{L^\infty(B_r(z))}$ is that even if $\norm{u^+}_{L^\infty(B_r(z))} > \norm{u^-}_{L^\infty(B_r(z))}$, we may have $\norm{u^+}_{L^\infty(B_r(z))} \ll r^{1-p/q} \norm{u^-}_{L^\infty(B_r(z))}^{p/q}$, or vice versa.
\end{rem} 

Next, we observe that each phase is a weak subsolution. 

\begin{lem}\label{lem:basic0}
Let $u\in W^{1,p\wedge q}(\Omega)$ be a local minimizer. Then $u^+ \in W^{1,p}(\Omega)$ is a weak $p$-subsolution and $u^- \in W^{1,q}(\Omega)$ is a weak $q$-subsolution.
\end{lem}

\begin{proof}
Since the argument is symmetric in between $u^+$ and $u^-$, we shall only prove that for $u^+$. Fix $\delta>0$, and consider a smooth approximation $\beta_\delta$ of the Heaviside function on the real line. That is, $\beta_\delta \in C^\infty((-\infty,\infty))$ such that
$$
\begin{cases}
\beta_\delta (t) = 0 & \text{for } t < 2^{-1}\delta,\\
\beta_\delta (t) > 0 & \text{for } t \geq 2^{-1}\delta,\\
\beta_\delta (t) = 1 & \text{for } t>\delta,\\
\beta_\delta' \geq 0 & \text{for } t\in(-\infty,\infty). 
\end{cases}
$$
Let $\eta\in C_0^\infty(\Omega)$ be nonnegative, and define $\phi = \eta \beta_\delta(u)$. Let $\e>0$ be any small number satisfying $2\e\norm{\eta}_{L^\infty(\Omega)} \leq \delta$. Clearly, $u - \e\phi \in W_{loc}^{1,1}$ and $\supp\phi \subset \Omega$, so the minimality of $u$ implies that $J( u - \e\phi,\supp \phi) \geq J(u,\supp\phi)$. Note that $\supp \phi \subset \supp\eta \cap \overline{\{u> \frac{\delta}{2}\}}$, since $\beta_\delta (u) > 0$ in $\{ u > \frac{\delta}{2}\}$. Also remark that $\{ u < \e\phi\}\subset\{u\leq 0\}$, since $\e \| \phi\|_{L^\infty(\Omega)} \leq \e \|\eta\|_{L^\infty(\Omega)} \leq \frac{\delta}{2}$ and $\phi = 0$ on $\{u < \frac{\delta}{2}\}$. Therefore, we have $(u-\e\phi)^+ = u-\e\phi$ and $(u-\e\phi)^- = 0$ in $\{u>\frac{\delta}{2}\}$. Then it follows from the minimality condition that 
$$
\begin{aligned}
0 &\leq \int_{\{u>\delta/2\}} \frac{|D(u-\e\phi)|^p - |Du|^p}{\e} \,dx \\
&= - \int_{\{u>\delta/2\}} |Du|^{p-2}\langle Du, D\phi\rangle\,dx + o(1)\\
&= - \int_{\{u>\delta/2\}} \beta_\delta(u)|Du|^{p-2}\langle Du, D\eta\rangle\,dx + o(1),
\end{aligned} 
$$  
where $o(1)$ is a term tending to zero as $\e\to 0$. Letting $\e\to 0$, then sending $\delta\to 0$, and utilising the convergence $\beta_\delta(u)\to \chi_{\{u>0\}}$, which holds a.e.\ on $\Omega$, we may deduce from the dominated convergence theorem that 
$$
\int_\Omega  |Du|^{p-2} \langle Du, D\eta\rangle\,dx \leq 0.
$$ 
Since $\eta$ is an arbitrary nonnegative function in $C_0^\infty(\Omega)$, we conclude that $u^+ \in W^{1,p}(\Omega)$ is a weak $p$-subsolution.
\end{proof}

Let us also state an approximation lemma that if $u^-$ is small, then $u^+$ is close to a $p$-harmonic function.

\begin{lem}\label{lem:basic1}
Let $u\in W^{1,p\wedge q}(B_4)$ be a local minimizer of $J$, and $v \in u^+ + W_0^{1,p}(B_1)$ be the $p$-harmonic function. Then 
\begin{equation}\label{eq:basic1}
\int_{B_1} |D (u^+ - v)|^p \,dx \leq c \begin{dcases}
c\| u^- \|_{L^q(B_4)}^q, & \text{if } p \geq 2, \\
c\norm{u^-}_{L^q(B_4)}^{\frac{qp}{2}} \norm{D u^+}_{L^p(B_1)}^{p(1-\frac{p}{2})} &\text{if } 1 < p <2,
\end{dcases}
\end{equation}
where $c$ depends only on $n$, $p$ and $q$.  
\end{lem}

\begin{proof}
Set $\kappa = \norm{u^-}_{L^q(B_1)}$. As $u^-$ being a weak $q$-subsolution, the Cacciopoli inequality yields
\begin{equation}\label{eq:basic1-Du--Lq}
\int_{B_2} |D u^-|^q\,dx \leq C\kappa^q.
\end{equation}
Let $v \in u^+ + W_0^{1,p}(B_1)$ be the $p$-harmonic function, and let $\eta\in C_0^\infty(B_4)$ be a cutoff function such that $0\leq \eta \leq 1$ in $B_4$, $\eta = 1$ on $B_1$, $|D \eta| \leq 4$ on $B_2$, and $\supp\eta \subset B_2$. Define an auxiliary function $\vp: \Omega \to \R$ by 
$$
\vp = 
\begin{dcases}
u &\text{in }B_4\setminus B_2\\
u^+ - (1-\eta) u^- &\text{in } B_2\setminus B_1,\\
v & \text{in }  B_1,
\end{dcases}
$$
Since $v - u^+ \in W_0^{1,p}(B_1)$, $\eta = 1$ on $\partial B_1$ and $\eta = 0$ on $\partial B_2$, $\vp \in W^{1,p\wedge q}(B_4)$. Since $|D \eta|\leq 4$ in $B_2$, we observe from \eqref{eq:basic1-Du--Lq} that
$$
\begin{aligned}
\int_{B_2\setminus B_1} |D ((1-\eta)u^-)|^q\,dx &\leq 2^q\int_{B_2\setminus B_1}\left(|D u^-|^q + (u^-)^q|D \eta|^q \right)dx  \leq C\kappa^q.
\end{aligned}
$$
Hence, it follows immediately that 
\begin{equation}\label{eq:basic1-vp-J}
\begin{aligned}
J(\vp, B_2) & \leq \int_{B_1} |Dv|^p \,dx  +   \int_{B_2\setminus B_1} |Du^+|^p \,dx  + C\kappa^q.
\end{aligned}
\end{equation} 

On the other hand, since $\supp(\vp - u) = \overline{B_2} \subset \Omega$, $J_p(u,B_2;a_p,a_q) \leq J_q(\vp,B_2;a_p,a_q)$, which along with \eqref{eq:basic1-vp-J} implies that 
$$
J(u,B_2) \leq \int_{B_1} |D v|^p \,dx  +   \int_{B_2\setminus B_1} |D u^+|^p \,dx  + C\kappa^q.
$$ 
However, by the minimality of $ \int_{B_1} |Dv|^p \,dx$ among all functions in $u^+ + W_0^{1,p}(B_2)$, we can proceed as  
\begin{equation}\label{eq:basic1-u-v-1}
0\leq \int_{B_1} (|Du^+|^p- |Dv|^p)\,dx \leq C\kappa^q. 
\end{equation}
The conclusion now follows from some manipulations with elementary vector inequalities, see \cite[Page 100]{DP}. 
\end{proof} 

Employing the decay of the minimizers of functionals with standard $p$-growth, we can also derive a similar estimate for the $p$-th energy of $u^+$. The estimate below becomes valuable when the size of $u^-$ is small in a large ball. We will encounter such a situation later in Section \ref{section:holder}. 

\begin{lem}\label{lem:basic2}
Let $u \in W^{1,p\wedge q}(B_2)$ be a local minimizer. Then 
\begin{equation}\label{eq:basic2}
\int_{B_r} |D u^+|^p\,dx \leq c \int_{B_2} \left( r^{n-p + \sigma p} |D u^+|^p + (u^-)^q \right)dx,\quad\forall r\in(0,1),
\end{equation} 
where $\sigma \in (0,1)$ depends only on $n$, $p$, and $c> 1$ may depend further on $q$. 
\end{lem}

\begin{proof}
Put $\kappa = \norm{u^-}_{L^q(B_2)}$. Choose $v \in u^+ + W_0^{1,p}(B_1)$ to be the $p$-harmonic function. By \eqref{eq:basic1-u-v-1}, the minimality of $\int_{B_1}|Dv|^p\,dx$ and \cite[Theorem 7.7]{Giu},
$$
\begin{aligned} 
\int_{B_r} |D u^+|^p\,dx &\leq \int_{B_r} |Dv|^p\,dx + C\kappa^q \\
&\leq cr^{n-p + \sigma p} \int_{B_1} |Dv|^p\,dx + C\kappa^q \\
&\leq cr^{n-p + \sigma p} \int_{B_1} |Du^+|^p\,dx + C\kappa^q,
\end{aligned} 
$$
for any $r\in(0,1)$, with $\sigma$ depending only on $n$ and $p$.  
\end{proof} 

To the rest of this section, we shall present some more advanced properties of the minimizers. These results will play important roles in Section \ref{section:holder} and Section \ref{section:nonflat}. 

Based on the previous two lemmas, we shall observe that the measure of the set where $u$ is not large can be made small, if both the $L^q$-norm of $u^-$ and the $L^p$-norm of $|D u^+|$ is sufficiently small compared to the $L^p$-norm of $u^+$, up to a correct scaling factor. 

\begin{lem}\label{lem:basic3}
Let $u\in W^{1,p\wedge q}(B_4)$ be a local minimizer such that
\begin{equation}\label{eq:basic3-assump}
\mean{B_4} (u^+)^p\,dx = 1,\quad \mean{B_4} ((u^-)^q + |D u^+|^p)\,dx \leq \kappa,
\end{equation}
for some $\kappa>0$. Then 
\begin{equation}\label{eq:basic3}
\frac{| \{u < \frac{1}{2} \} \cap B_1|}{|B_1|} \leq c\kappa,
\end{equation}
where $c$ depends only on $n$, $p$, $q$ and $\Lambda$. 
\end{lem}

\begin{proof}
Let $v\in u^+ + W_0^{1,p}(B_2)$ be the minimizer of $J_+ \equiv J (\cdot,a_p,a_p)$. By Lemma \ref{lem:basic2} and the assumption \eqref{eq:basic3-assump}, we can deduce 
$$
\int_{B_2} |D(u^+ - v)|^p\,dx \leq c \kappa,
$$
for any $p >1$. Since $u^+ - v \in W_0^{1,p}(B_2)$, the Poincar\'e inequality implies
\begin{equation}\label{eq:u+-v-Lp}
\int_{B_2} |(u^+ - v)|^p\,dx \leq c \kappa. 
\end{equation}  
By \eqref{eq:basic3-assump}, we may proceed with a compactness argument based on the Sobolev embedding theory to choose $\bar\kappa$ sufficiently small (depending only on $n$ and $p$) such that
\begin{equation}\label{eq:u+-Lp}
\mean{B_1} (u^+)^p\,dx \geq \frac{7^p}{8^p}.
\end{equation}
Taking $\bar\kappa$ smaller if necessary, \eqref{eq:u+-v-Lp} and \eqref{eq:u+-Lp} ensures that 
\begin{equation}\label{eq:v-sup}
\sup_{B_1}v \geq\frac{3}{4}. 
\end{equation}
However, we can deduce from the interior H\"older estimate \cite[Theorem 2.9]{L} and the minimizing property of $v$, 
\begin{equation}\label{eq:v-Ca}
\osc_{B_1} v \leq  c\kappa^{\frac{1}{p}}. 
\end{equation}
Therefore, if $\bar\kappa$ is even smaller, then by \eqref{eq:v-sup} and \eqref{eq:v-Ca}, 
\begin{equation}\label{eq:v-inf}
\inf_{B_1} v \geq \frac{5}{8}.  
\end{equation}
Thus, we may derive from \eqref{eq:u+-v-Lp} and \eqref{eq:v-inf} that
$$
\begin{aligned}
\frac{|\{ u \leq \frac{1}{2}\}\cap B_1|}{|B_1|}  &\leq 8^p \int_{\{ u < \frac{1}{2}\}\cap B_1} \left(\frac{5}{8} - u^+\right)^p dx \\
&\leq 8^p \int_{\{ u \leq \frac{1}{2}\}\cap B_1} (v - u^+)^p \, dx\\
&\leq 8^p \int_{B_1} | v- u^+|^p\,dx\\
&\leq c\kappa, 
\end{aligned}
$$ 
which finishes the proof.
 \end{proof}

The next lemma asserts that the $L^q$-norm of $u^-$ in the interior can be controlled by that of $u^+$ over a larger ball, provided that $\{u \leq 0\}$ occupies positive measure in the interior. It is noteworthy that the control is independent of the $L^\infty$-norm of $u^-$ over the larger ball. 
\begin{lem}\label{lem:holder3}
Let $u\in W^{1,p\wedge q}(B_4)$ be a local minimizer such that 
\begin{equation}\label{eq:u-msr-var-re}
\frac{|\{ u\leq 0\}\cap B_1|}{|B_1|} \leq \gamma, 
\end{equation}
for some $\gamma \in (0,1)$. Then 
\begin{equation}\label{eq:holder3-re}
\| u^- \|_{L^q(B_1)} \leq \frac{c}{(1-\gamma)^{\frac{2}{q}\vee 1}} \| u^+\|_{L^p(B_4)}^{\frac{p}{q}}. 
\end{equation}
\end{lem}

\begin{proof}
Throughout the proof, $c$ will be a constant depending at most on $n$, $p$, $q$ and $\Lambda$, and it may be different at each occurrence. We may assume without loss of generality that $\mean{B_4} (u^+)^p \,dx = 1$, since the general case can be deduce by considering $w = \frac{u^+}{S_+} - \frac{u^-}{S_+^{p/q}}$ with $S_+ = (\mean{B_4}(u^+)^p\,dx)^{1/p}$. Write 
$$
A := \left( \mean{B_1} (u^-)^q\,dx \right)^{\frac{1}{q}}. 
$$
Assume that $A > 0$, since otherwise the proof becomes trivial.

We shall divide the proof into two cases, (i) $2\leq q<\infty$ and (ii) $1<q<2$, since the gradient difference between $u^-$ and its $q$-harmonic replacement satisfies different inequalities (see Lemma \ref{lem:basic1}) depending on the range of $q$. The idea of the proof is the same, but the argument for the case $1<q<2$ is technically more involved than the case  $2\leq q<\infty$.

\begin{case} 
$2\leq q<\infty$.
\end{case}

Let $v \in u^- + W_0^{1,q}(B_2)$ be the $q$-harmonic function. As $u^-\in W^{1,q}(B_4)$ being a weak $a_q$-subsolution (Lemma \ref{lem:basic0}), it follows from the comparison principle that $v \geq u^-$ a.e.\ on $B_2$. Thus, one can deduce from the Harnack inequality that
$$
c\inf_{B_1} v \geq \sup_{B_1} v \geq A,
$$
where $c>1$ depends only on $n$ and $q$. Then we obtain 
\begin{equation}\label{eq:holder3-re-1}
\begin{aligned} 
\left| \left\{ u^- < \frac{A}{2c} \right\} \cap B_1 \right| &\leq \left(\frac{2c}{A}\right)^q \int_{B_1 \cap \{u^- < \frac{A}{2c}\}} \left(\frac{A}{c} - u^-\right)^q\,dx\\
&\leq \left(\frac{2c}{A}\right)^q \int_{B_1\cap \{u^- < \frac{A}{2c}\}} ( v- u^- )^q\,dx \\
&\leq \left(\frac{2c}{A}\right)^q \int_{B_2} ( v - u^- )^q \,dx. 
\end{aligned}
\end{equation}
Thanks to \eqref{eq:u-msr-var-re}, we have $|\{ u \geq 0\}\cap B_1| \geq (1-\gamma)|B_1|$. Hence, it follows from \eqref{eq:holder3-re-1}, $\{ u \geq 0\}\subset \{ u^- < \frac{A}{2c}\}$ and the Poincar\'e inequality that 
\begin{equation}\label{eq:a-hol3}
\begin{aligned}
A^q \leq \frac{c}{1-\gamma} \int_{B_2} |D (v-u^-)|^q\,dx.
\end{aligned}
\end{equation}  

Note that $(-u)$ is a local minimizer to the functional $J_{q,p}$, according to Lemma \ref{lem:basic}. Since $q\geq 2$, we may apply Lemma \ref{lem:basic1} to $(-u)$ and deduce from the assumption $\mean{B_4} (u^+)^p\,dx = 1$ that 
$$
\int_{B_2} |D (u^- - v)|^q\,dx \leq c.
$$
Inserting this inequality to \eqref{eq:a-hol3} yields 
$$
A^q \leq \frac{c}{1-\gamma}
$$
which proves \eqref{eq:holder3-re} for case $2\leq q < \infty$. 

\begin{case}
$1<q<2$. 
\end{case}

As for the case $1<q<2$, let us fix $1< s< r<4$, and write 
$$
A_s := \left(\mean{B_s} (u^-)^q\,dx \right)^{\frac{1}{q}},\quad t:= \frac{r+s}{2}, 
$$ 
and again assume without loss of generality that $A_s>0$. 

Let $v_t \in u^- + W_0^{1,q}(B_t)$ be the $q$-harmonic function. Following the above argument for the case $2\leq q<\infty$, we also obtain that $v_t \geq u^-$ a.e.\ on $B_t$, and that
\begin{equation}\label{eq:harn-q}
 \frac{c}{(r-s)^\sigma} \inf_{B_s} v_t \geq \sup_{B_s} v_t \geq A_s,
\end{equation}
where $\sigma>0$ is a constant depending only on $n$ and $q$. We may proceed analogously as with the computation in \eqref{eq:holder3-re-1} and obtain that 
\begin{equation}\label{eq:holder3-q-1}
\left| \left\{ u^- < \frac{(r-s)^\sigma A_s}{2c} \right\} \cap B_1 \right| \leq \left(\frac{2c}{(r-s)^\sigma A_s}\right)^q \int_{B_t} (v- u^-)^q\,dx. 
\end{equation} 
Thus, as in \eqref{eq:a-hol3}, we may also deduce from \eqref{eq:u-msr-var-re} and the Poincar\'e inequality that 
\begin{equation}\label{eq:holder3-q-2}
A_s^q \leq \frac{c}{(1-\gamma)(r-s)^{\sigma q}} \int_{B_t} |D (v-u^-)|^q\,dx. 
\end{equation} 

Now that we assume $1<q<2$ and $s<t<r$, we may apply Lemma \ref{lem:basic1} and the Cacciopoli inequality to $(-u)$ in the rescaled form and derive from $\mean{B_4}(u^+)^p\,dx = 1$ that
$$
\begin{aligned}
\int_{B_t} |D (u^- - v)|^q\,dx &\leq \frac{c }{(r-s)^p} \norm{D u^-}_{L^q(B_t)}^{q(1-\frac{q}{2})}\leq \frac{c A_r^{q(1-\frac{q}{2})}}{(r-s)^{p + q(1-\frac{q}{2})}},
\end{aligned}
$$
where we used the Caccioppoli inequality to derive the rightmost side. Using this inequality to substitute the rightmost integral in \eqref{eq:holder3-q-2}, we arrive at
\begin{equation}\label{eq:holder3-q-3}
\begin{aligned}
A_s^q &\leq \frac{c}{(1-\gamma)(r-s)^{(\sigma + 1 - \frac{q}{2})q + p}} A_r^{q(1-\frac{q}{2})} \\
& \leq \frac{c}{(1-\gamma)^{\frac{2}{q}}(r-s)^{2(\sigma + 1 - \frac{q}{2}) + \frac{2p}{q}}} + \frac{1}{2} A_r^q,
\end{aligned}
\end{equation} 
where the last inequality follows from Young's inequality.

As $|B_s|A_s^q = \| u^-\|_{L^q(B_s)}^q$ being nondecreasing in $s$, we may now invoke a standard iteration lemma  \cite[Lemma 6.1]{Giu} to deduce 
$$
A^q \leq  \frac{c}{(1-\gamma)^{\frac{2}{q}}},
$$
where $A = A_1$. This finishes the proof. 
\end{proof}

As a corollary to Lemma \ref{lem:holder3}, we observe that the $L^q$-norm of $u^-$ in an interior ball can be bounded by a universal constant multiple of the $L^p$-norm of $u^+$ over the same ball, provided that $u^+$ satisfies a doubling condition. Here we present the statement in a more general context for the future reference.

\begin{cor}\label{cor:holder3}
Under the assumption of Lemma \ref{lem:holder3}, suppose further that 
$$
\norm{u^+}_{L^p(B_4)} \leq \frac{1}{\beta} \norm{u^+}_{L^p(B_1)},
$$
for some $\beta\in (0,1)$. Then  
\begin{equation}\label{eq:holder3}
\frac{\norm{u^-}_{L^\infty(B_1)}}{\norm{u^+}_{L^\infty(B_1)}^{\frac{p}{q}}}\leq \frac{C}{\beta^{\frac{p}{q}}(1-\gamma)^{\frac{2}{q}\vee 1}} ,
\end{equation} 
where $C$ and $\mu$ are positive constants depending only on $n$, $p$ and $q$. 
\end{cor}


\section{H\"older Continuity  (Proof of Theorem \ref{thm:holder})}\label{section:holder} 

This section is devoted to interior H\"older regularity for local minimizers. By the Sobolev theory, the theorem is trivial if $\min\{p,q\}>n$. Nevertheless, we shall {\it not} restrict the range of $p$ and $q$, unless stated otherwise. 

Let $\Gamma$ be the set of all ``vanishing points'' of $u$ in the Lebesgue sense: 
\begin{equation}\label{eq:gam}
\Gamma = \left\{ z \in \Omega: \lim_{r\to 0} \mean{B_r(z)} u\,dx = 0\right\}.
\end{equation} 
We consider the vanishing points in the Lebesgue sense because we do not have continuity of minimizers yet. 

Our primary goal is to prove the H\"older growth of the size of $u^+$ ($u^-$) around each ``nonpositive'' (resp. nonnegative) point of $u$. As the argument being symmetric, we shall only present the assertion for $u^+$. The H\"older growth will be obtain by an iteration of the following proposition.

\begin{prop}\label{prop:holder}
There exists a constant $0<\eta<\frac{1}{2}$, depending only on $n$, $p$ and $q$, such that if $u\in W^{1,p\wedge q}(B_4)$ is a local minimizer such that 
$$
\mean{B_4} (u^+)^p\,dx = 1,\quad \limsup_{r\to 0}\mean{B_r} u\,dx \leq 0,
$$
then 
$$
\mean{B_1} (u^+)^p \,dx \leq 1-\eta.
$$
\end{prop}

Observe that as $u^+$ being a weak $a_p$-subsolution, the conclusion of the above proposition will follow immediately from the De Giorgi lemma \cite{Giu}, once we prove that the set $\{u<\frac{1}{2}\}\cap B_1$ contains some universal amount of measure in $B_1$. For this purpose, we shall consider its contraposition. 

\begin{prop}\label{prop:holder-re}
There exist $\tau > 0$ and $\rho > 0$, depending only on $n$, $p$ and $q$, such that if $u\in W^{1,p\wedge q}(B_4)$ is a local minimizer satisfying  
$$
\mean{B_4} (u^+)^p\,dx =1,\quad \frac{|\{ u < \frac{1}{2}\} \cap B_2|}{|B_2|} \leq \tau,
$$
then $u > 0$ a.e.\ in $B_\rho$.
\end{prop}

This proposition will be proved by the following key lemma.

\begin{lem}\label{lem:pos}
There exists $\e>0$, depending only on $n$, $p$ and $q$, such that if $u\in W^{1,p\wedge q}(B_4)$ is a local minimizer satisfying
$$
\mean{B_4} (u^+)^p\,dx = 1,\quad\mean{B_4} ((u^-)^q + |Du^+|^p)\,dx \leq\kappa,
$$
then $u> 0$ a.e.\ in $B_1$. 
\end{lem} 

\begin{proof}
Let us consider the case $q < n$. The argument for the case $q \geq n$ is almost the same except for the part where we use the Sobolev embedding. Hence, we shall mention the change to be made for the latter case, at the end of the proof, and leave out the detail to the reader. 

Let $\sigma$ be as in Lemma \ref{lem:basic2}, and choose $s\in(0,\frac{1}{8})$ such that 
\begin{equation}\label{eq:s}
cs^{\frac{\sigma p}{2}} \leq \frac{1}{2^{p+2}},\quad s^{\frac{n^2}{q} - n} \leq 2^{2n-1}, 
\end{equation}
for some large constant $c>1$ depending only on $n$, $p$ and $q$. We shall take $\e$ sufficiently small such that 
\begin{equation}\label{eq:t-k}
c\e^{\frac{q}{n}} \leq \frac{s^{n + (q- (1- \frac{\sigma}{2} )p)^+}}{2^{p+2}}.
\end{equation}

For each $\ell\in \N\cup\{0\}$, define $u_\ell : B_4\to \R$ by 
\begin{equation}\label{eq:ul}
u_\ell (x) = \frac{u^+(s^\ell x)}{(4s^\ell)^{-\frac{n}{p}}\| u^+\|_{L^p(B_{4s^\ell})}} - \frac{u^-(s^\ell x)}{4^{-\frac{n}{q}} s^{\ell(1-\frac{p}{q} - \frac{n}{q})}\| u^+\|_{L^p(B_{4s^\ell})}^{\frac{p}{q}}}.
\end{equation}
Then $u_\ell \in W^{1,p\wedge q}(B_4)$ is a cubical $(1+\kappa)$-minimum of the functional $J$. Suppose that for some $\ell \in \N\cup\{0\}$, the hypothesis of this lemma holds with $u_\ell$ in place of $u$, i.e., 
\begin{equation}\label{eq:ul-Lp-Lq}
\mean{B_4} (u_\ell^+)^p \,dx = 1,\quad \mean{B_4} ((u_\ell^-)^q + |Du_\ell^+|^p) \,dx \leq s^{\frac{\ell\sigma p}{2}}\kappa. 
\end{equation}
Note that this hypothesis holds for $\ell = 0$ with $u_0 = u$. 

First, we observe from Lemma \ref{lem:basic3} that 
\begin{equation}\label{eq:ul-msr}
\frac{| \{ u_\ell \leq \frac{1}{2} \} \cap B_1|}{|B_1|} \leq cs^{\frac{\ell\sigma p}{2}}\kappa.  
\end{equation} 
By \eqref{eq:s} and \eqref{eq:t-k}, we have $c\kappa\leq s^{\frac{n^2}{q}} \leq \frac{1}{2}|B_{4s}|$. Therefore, by \eqref{eq:ul-msr},
\begin{equation}\label{eq:ul+-Lp}
\mean{B_{4s}} (u_\ell^+)^p\,dx \geq \frac{|\{ u_\ell \geq \frac{1}{2}\}\cap B_{4s}|}{2^p |B_{4s}|} \geq \frac{1}{2^{p+1}}. 
\end{equation}
By the Cacciopoli inequality (as $u_\ell^-$ being a weak $a_q$-subsolution), we have 
\begin{equation}\label{eq:Dul--Lq}
\int_{B_1} |Du_\ell^-|^q \,dx \leq cs^{\frac{\ell\sigma p}{2}}\kappa. 
\end{equation}
Therefore, applying H\"older inequality, the Sobolev inequality (of $L^{\frac{nq}{n-q}} \hookrightarrow L^q$) and the Poincar\'e inequality, we may proceed by means of together with \eqref{eq:ul-msr} and \eqref{eq:Dul--Lq} as 
\begin{equation}\label{eq:ul--Lq-re}
\begin{aligned}
\int_{B_1} (u_\ell^-)^q \,dx &\leq c (s^{\frac{\ell\sigma p}{2}}\kappa)^{\frac{q}{n}} \left( \int_{B_1} (u_\ell^-)^{\frac{nq}{n-q}}\,dx\right)^{1-\frac{q}{n}} \\
&\leq c (s^{\frac{\ell\sigma p}{2}}\kappa)^{\frac{q}{n}} \int_{B_1} ( |Du_\ell^-|^q + (u_\ell^-)^q)\,dx \\
& \leq c(s^{\frac{\ell\sigma p}{2}}\kappa)^{\frac{q}{n}} \int_{B_1} |Du_\ell^-|^q\,dx \\
&  \leq c(s^{\frac{\ell\sigma p}{2}}\kappa)^{1+\frac{q}{n}} ;
\end{aligned}
\end{equation}
here the Poincar\'e inequality was used in the derivation of the third inequality, and was available because $|\{ u_\ell \geq 0 \}\cap B_1| > 1 - cs^{\ell \sigma p}\kappa> \frac{1}{2}$, thanks to \eqref{eq:ul-msr}, \eqref{eq:t-k} and \eqref{eq:s}. As a byproduct of \eqref{eq:t-k} and \eqref{eq:ul+-Lp}, 
\begin{equation}\label{eq:ul--Lq}
\begin{aligned}
\mean{B_{4s}} (u_\ell^-)^q \,dx  &\leq c  \e^{1 + \frac{q}{n}} s^{-n + \ell (1+\frac{q}{n})\frac{\sigma p}{2} } \\
& \leq \frac{\kappa ^{q- (1 -(\ell + 1) \frac{\sigma }{2})p}}{2^{p+2}} \\
&\leq \frac{\kappa ^{q- (1  - (\ell + 1) \frac{\sigma }{2})p}}{2}  \mean{B_{4s}} (u_\ell ^+)^p\,dx. 
\end{aligned}
\end{equation}
Besides, by Lemma \ref{lem:basic2}, \eqref{eq:ul-Lp-Lq}, \eqref{eq:ul--Lq-re}, \eqref{eq:s}, \eqref{eq:t-k} and \eqref{eq:ul+-Lp},
\begin{equation}\label{eq:Dul+-Lp}
\begin{aligned}
\mean{B_{4s}} |Du_\ell^+|^p \,dx & \leq c s^{-n + \frac{\ell\sigma p}{2}} \kappa( (s^{n - p + \sigma p} + \kappa) + (s^{\ell\sigma p}\kappa)^{\frac{q}{n}}) \\
&\leq c \kappa ^{-(1 - \sigma  - \frac{\ell \sigma}{2}) p} \\
&\leq \frac{ \kappa ^{-(1  - (\ell +1)\frac{\sigma}{2}) p}}{2^{p+2}} \\
&\leq \frac{\kappa ^{-(1 - (\ell +1)\frac{\sigma}{2}) p} }{2}\mean{B_{4s}} (u_\ell^+)^p\,dx.
\end{aligned}
\end{equation}

Define $u_{\ell+1} : B_4 \to \R$ as in \eqref{eq:ul} with $\ell+1$ in place of $\ell$. One may verify directly from \eqref{eq:ul--Lq} and \eqref{eq:Dul+-Lp} that \eqref{eq:ul-Lp-Lq} holds with $\ell$ replaced by $\ell+1$. By the induction principle, we observe that \eqref{eq:ul-Lp-Lq} holds for every number $\ell\in\N\cup\{0\}$. In particular, rephrasing \eqref{eq:ul-Lp-Lq} in terms of $u$ yields (after some obvious manipulations of the constant) that 
\begin{equation}\label{eq:u-Lpq-gr}
\mean{B_r} |Du^+|^p dx \leq c\kappa r^{-(1-\frac{\sigma}{2})p} ,\quad\forall r\in(0,1),
\end{equation}
Here we also used the local $L^\infty$-estimate (Lemma \ref{lem:basic0}) for $u^+$, which along with $\mean{B_4} (u^+)^p\,dx = 1$ implies $\| u^+\|_{L^\infty(B_2)} \leq c$. 

Taking $\kappa$ even smaller if necessary, we may repeat the above argument around any point $z\in B_1$, and obtain
$$
\mean{B_r(z)} |Du^+|^p \,dx \leq c\kappa r^{-(1-\frac{\sigma}{2})p},\quad\forall r\in(0,1),\,\forall z\in B_1,
$$
possibly with a larger constant $c$. Therefore, by Morrey's lemma, we deduce that $u^+ \in C^{0,\frac{\sigma}{2}}(B_1)$ and
\begin{equation}\label{eq:u+-Ca-re}
[u^+]_{C^{0,\frac{\sigma}{2}}(B_1)} \leq c\kappa^{\frac{1}{p}}. 
\end{equation}
Finally, by \eqref{eq:ul-msr} with $\ell = 0$, we have $|\{ u \leq \frac{1}{2}\}\cap B_1| \leq c\e$. Hence. with $c\kappa\leq 2^{-2n -1}$, we have $|\{ u > \frac{1}{2}\}\cap B_1| > 0$, which now implies via \eqref{eq:u+-Ca-re} that 
$$
\inf_{B_1} u^+ \geq \frac{1}{2} - c \kappa^{\frac{1}{p}} > 0,
$$
provided that we choose $\kappa$ even smaller. This finishes the proof for the case $q < n$. 

As for the case $q = n$, we only need to replace $\frac{q}{n}$ at each point of its appearance (which is as a matter of fact $1- \frac{q}{q_*}$, with $q_* = \frac{nq}{n-q}$ being the Sobolev conjugate when $q< n$) with some constant $1-\nu$, with $\nu$ sufficiently small. Especially, one needs to modify \eqref{eq:s} and \eqref{eq:t-k} in this way. Then the estimate \eqref{eq:ul--Lq} would also hold, by utilizing a modified Sobolev embedding that $\|u_\ell^-\|_{L^{\frac{q}{\nu}}(B_{1/2})} \leq c_\nu \|u_\ell^-\|_{W^{1,q}(B_{1/2})} \leq c_\nu$. The estimate \eqref{eq:Dul+-Lp} should hold in the same way, once we assume in \eqref{eq:s} that $s^{\frac{n}{1-\nu}} \leq \frac{1}{2} s^n$; as a matter of fact, here is where the size of $\nu$ can be deduced. The rest of the argument can be repeated verbatim, whence we omit the detail. 

As for the case $q > n$, we remove the condition $s^{\frac{n^2}{q} - n} \leq \frac{1}{2}$ in \eqref{eq:s} and replace $\kappa^{\frac{q}{n}}$ in \eqref{eq:t-k} with $\kappa$. Moreover, we replace the $L^q$-norm of $u_\ell^-$ by $L^\infty$-norm at each occurrence; note by the Sobolev embedding $u_\ell^- \in C^{0,1-\frac{n}{q}}(B_2)$ in this case. The argument can then be repeated in the same way, without any major modification, so we skip the detail. 
\end{proof} 

To fall into the setting of Lemma \ref{lem:pos} from that of Proposition \ref{prop:holder-re}, we use the comparability of two phases, Corollary \ref{lem:holder3}. 

\begin{proof}[Proof of Proposition \ref{prop:holder-re}]
By the assumption $|\{ u \leq \frac{1}{2}\}\cap B_2| < \tau |B_2|$, choosing $\tau$ sufficiently small, depending only on $n$, we may deduce analogously as with \eqref{eq:ul-msr} that 
$$
\mean{B_1} (u^+)^p \,dx \geq \frac{1}{4^p}. 
$$
As we assume $\mean{B_4} (u^+)^p\,dx = 1$, we may use Corollary \ref{cor:holder3} (with $\beta = \frac{1}{4}$ and $\gamma \leq 4^n\tau < \frac{1}{2}$) after suitable scaling that 
\begin{equation}\label{eq:uB-1}
\norm{u^-}_{L^q(B_1)} \leq c\norm{u^+}_{L^p(B_1)}^{\frac{p}{q}},
\end{equation} 
where $c$ depends only on $n$, $p$ and $q$. To the rest of the proof, $c$ may vary at each appearance. 

Substituting \eqref{eq:ul-msr} with the our assumption on the measure of $\{ u \leq\frac{1}{2}\}$ and using \eqref{eq:uB-1}, we may follow the derivation of \eqref{eq:ul+-Lp}, \eqref{eq:ul--Lq} and \eqref{eq:Dul+-Lp} and find a sufficiently small $\rho \in (0,\frac{1}{8})$, yet depending only on $n$, $p$ and $q$, such that 
$$
\mean{B_{4\rho}} \left(\frac{(u^-)^q}{\rho^q} + |D u^+|^p\right)dx \leq c\tau \rho^{q-p} \mean{B_{4\rho}} (u^+)^p\,dx.
$$
Thus, we may rescale our minimizer $u$ as 
$$
u_\rho (x) = \frac{u^+(\rho x)}{|B_{4\rho}|^{-\frac{1}{p}}\|u^+\|_{L^p(B_{4\rho})}} - \frac{u^-(\rho x)}{\rho^{1-\frac{p}{q} - \frac{n}{q}}|B_{4\rho}|^{-\frac{1}{q}}\|u^+\|_{L^p(B_{4\rho})}},
$$
we see that $u_\rho \in W^{1,p\wedge q}(B_4)$ is a minimizer of the functional $J$, such that
$$
\mean{B_4} (u_\rho^+)^p\,dx = 1,\quad \mean{B_4} ( (u_\rho^-)^q + |D u_\rho^+|^p) \,dx \leq c\tau. 
$$
Choosing $c\tau \leq \kappa$, with $\kappa$ as in Lemma \ref{lem:pos}, we obtain
$$
u_\rho > 0 \quad\text{a.e.\ in }B_1.
$$
Rescaling back, we obtain $u > 0$ a.e.\ in $B_\rho$, as desired.
\end{proof}

As mentioned in the beginning of this section, Proposition \ref{prop:holder} follows immediately from Proposition \ref{prop:holder-re} as a contraposition, along with the De Giorgi lemma \cite{Giu}. As a corollary, we obtain the H\"older growth at each ``non-positive'' point. 

\begin{cor}\label{cor:gam}
Let $u \in W^{1,p\wedge q}(\Omega)$ be a local minimizer such that 
$$
\limsup_{r\to 0} \mean{B_r(z)} u\,dx \leq 0. 
$$
for some $z\in \Omega$. Then 
\begin{equation}\label{eq:gam-decay}
\norm{u^+}_{L^\infty(B_r(z))} \leq \frac{Cr^\alpha}{(\dist(z,\partial\Omega))^{\alpha + \frac{n}{p}}} \norm{u^+}_{L^p(\Omega)},
\end{equation}
where $\alpha\in(0,1)$ and $C>0$ depend only on $n$, $p$ and $q$. 
\end{cor} 

We remark that the iteration with Proposition \ref{prop:holder} works as the proposition only requires the normalization of $L^p$-norm for $u^+$, and leaves $u^-$ free. Otherwise, the iteration may fail, as each phase scales differently; see Remark \ref{rem:scale}.

Now define 
\begin{equation}\label{eq:D+}
D^+ :=  \left\{z\in\Omega: \liminf_{r\to 0}  \mean{B_r(z)} u\,dx >0 \right\}, 
\end{equation}
and 
\begin{equation}\label{eq:D-}
D^- := \left \{z\in\Omega: \limsup_{r\to 0}  \mean {B_r(z)} u \,dx <0\right\}.
\end{equation}
Note that $D^+$ and $D^-$ are the sets of positive and respectively negative values of $u$ in the Lebesgue sense. In the following lemma, we assert that $D^+$, $D^-$ are open, $\Gamma$ is closed, and the given domain $\Omega$ is a disjoint union of $D^+$, $D^-$ and $\Gamma$. We shall skip the proof, as it only involves some elementary real analysis, and basic properties of $a_p$-harmonic functions. 

\begin{lem}\label{lem:gam-top}
Let $u \in W^{1,p\wedge q}(\Omega)$ be a local minimizer, and let $\Gamma$, $D^+$ and $D^-$ be defined as in \eqref{eq:gam}, \eqref{eq:D+} and \eqref{eq:D-}. Then $D^+$ and $D^-$ are open, and $\Gamma = \partial D^+\cap \Omega = \partial D^- \cap \Omega$. In particular,  $\Omega = D^+ \cup D^- \cup \Gamma$. 
\end{lem} 

\bigskip
\noindent
{\bf Proof of Theorem \ref{thm:holder}:}

We are now ready to prove the interior H\"older regularity of local minimizers. 

To simplify the exposition, let us assume that $\Omega = B_8$, $K = \overline{B_4}$ and $u \in W^{1,p\wedge q}(B_8)$ is a local minimizer such that $0\in\Gamma$, where $\Gamma$ is as in \eqref{eq:gam}. We shall only prove the estimate \eqref{eq:Ca} for $u^+$, as the argument for $u^-$ is symmetric. After scaling, one may also assume that 
$$
\mean{B_8} (u^+)^p\,dx = 1. 
$$

Let $x,y\in B_1\cap D^+$ be Lebesgue points of $u$ satisfying $|x-y| \leq  \min\{d_x,d_y\}$, where $d_x = \frac{1}{4}\dist(x,\Gamma)$ and $d_y = \frac{1}{4} \dist(y,\Gamma)$. Assume without loss of generality that $d_x \geq d_y$. Note that $x\in \overline{B_{d_y}(y)}\subset B_{4d_y}(y)\subset\Omega$. Also it follows from Corollary \ref{cor:gam} that 
$$
\norm{u^+}_{L^\infty(B_{2d_y}(y))} \leq cd_y^\sigma. 
$$
Consider an auxiliary function 
$$
v(\xi) := \frac{u^+(y + 2d_y \xi)}{c d_y^\sigma}.
$$
Observe that $v\in W^{1,p\wedge q}(B_4)$ is a positive $p$-harmonic function such that $\norm{v}_{L^\infty(B_4)} \leq 1$. Therefore, it follows from the interior H\"older estimate for $a_p$-harmonic functions that $v \in C^{0,\sigma}(\overline{B_1})$ and $[ v ]_{C^{0,\sigma}(\overline{B_1})} \leq c$, so 
$$
|u^+(x) - u^+(y)| \leq c |x-y|^\sigma,
$$ 
where $c$ and $\sigma$ depends only on $n$, $p$, $q$ and $\Lambda$.  

As for the case where $x,y\in B_1\cap D^+$ with $|x-y| > \min\{d_x,d_y\}$, the proof for the above inequality is more straightforward, and we shall skip the detail. This finish the proof, as $u^+ = 0$ a.e.\ in $B_1\setminus D^+$.


\section{Regularity of Flat Free Boundaries (Proof of Theorem \ref{thm:flat}) }\label{section:flat}

This section is devoted to the regularity of flat parts of the free boundary of local minimizers. Here and thereafter, given a local minimizer $u\in W^{1,p\wedge q}(\Omega)$, we shall call $\partial\{u > 0\}\cap\Omega$ the free boundary of $u$. Note that the free boundary of a local minimizer is well-defined, due to its (H\"older) continuity (Theorem \ref{thm:holder}).
The main result of this section is Theorem \ref{thm:flat}. The proof for this theorem will be divided into two steps: (i) local minimizers are viscosity solutions (in the sense below), and (ii) flat free boundaries of the viscosity solutions are locally $C^{1,\sigma}$-graphs. 

Let us first define the notion of viscosity solutions. 

\begin{defn}\label{defn:visc}
Given a domain $\Omega\subset\R^n$, we shall call $u\in C(\Omega)$ a viscosity solution, to the $(p,q)$-harmonic free boundary problem, if $u^+ \in W^{1,p}(\Omega)$ is $p$-harmonic in $\{u>0\}$, and $u^- \in W^{1,q}(\Omega)$ is $q$-harmonic in $\mbox{int}(\{ u \leq 0\} )$,\footnote{Since $u^-\in C(\Omega)\cap W^{1,q}(\Omega)$ is a $q$-harmonic in $\mbox{int}(\{ u \leq 0\})$, it follows from the strong minimum principle that $\mbox{int}(\{ u \leq 0\}) = \{ u < 0\}$, hence $\partial \{ u > 0\}\cap\Omega = \partial \{ u < 0\}\cap\Omega$.} and the following free boundary condition is satisfied for every $x_0\in\partial\{u>0\}\cap\Omega$ and $B\subset \{u>0\}\cup \mbox{int}(\{u \leq 0\})$ such that $\partial B \cap \partial\{u>0\} = \{x_0\}$. 
\begin{enumerate}[(a)] 
\item If $B\subset \{u>0\}$ with $\nu$ the unit inward normal to $\partial B$ at $x_0$, then there exist some $\alpha,\beta>0$, satisfying $(p-1)\alpha^p = (q-1)\beta^q$, such that 
\begin{align*}
u^+(x) &\geq \alpha \inn{x-x_0}{\nu}^+ + o(|x-x_0|)\quad\text{for }x\in B,\\
u^-(x) &\leq \beta\inn{x-x_0}{\nu}^- + o(|x-x_0|)\quad\text{for }x\in B^c
\end{align*}
as $x\to x_0$, with equality along every non-tangential domain in both cases. 
\item If $B \subset \mbox{int}(\{u \leq 0\})$ with $\nu$ the unit outward normal to $\partial B$ at $x_0$, then there exist some $\alpha,\beta>0$, satisfying $(p-1)\alpha^p = (q-1)\beta^q$, such that 
\begin{align*}
u^-(x) &\geq \beta \inn{x-x_0}{\nu}^- + o(|x-x_0|)\quad\text{for }x\in B,\\
u^+(x) &\leq \alpha \inn{x-x_0}{\nu}^+ + o(|x-x_0|)\quad\text{for } x\in B^c,
\end{align*}
as $x\to x_0$, with equality along every non-tangential domain in both cases.
\end{enumerate}
\end{defn} 

Our first main step is to prove the following proposition. 

\begin{prop}\label{prop:visc}
A local minimizer $u \in W^{1,p\wedge q}(\Omega)$ is a viscosity solution.
\end{prop} 

Once the proposition is proved, we shall prove the following proposition, which follows similar ideas as that of \cite[Theorem 1]{LN1} (that Lipschitz free boundaries are $C^{1,\sigma}$) and \cite[Theorem 1]{LN2} (that flat free boundaries are Lipschitz). 

\begin{prop}[Essentially due to {\cite{LN1,LN2}}]\label{prop:flat-fb}
Let $u \in C(\Omega)$ be a viscosity solution, and let $z\in\partial\{u>0\}\cap\Omega$ be arbitrary with $0<r<\dist(z,\partial\Omega)$. Then there exist some real numbers $\sigma\in(0,1)$, $\e>0$ and $\rho>0$, depending at most on $n$, $p$ and $q$, such that if $\partial\{u>0\}\cap B_r(z) \subset \{x\in\R^n: |\inn{x-z}{\nu}| \leq r\e\}$ for some direction $\nu\in\partial B_1$, then $\partial\{u>0\}\cap B_{\rho r}(z)$ is a $C^{1,\sigma}$-graph in direction $\nu$. 
\end{prop}

Note that Theorem \ref{thm:flat} would follow as a corollary from the above two propositions. 

Let us begin with the weak formulation of the free boundary of a local minimizer.

\begin{lem}\label{lem:euler}
Let $u \in W^{1,p\wedge q}(\Omega)$ be a local minimizer. Then $u^+ \in W^{1,p}(\Omega)$ is $p$-harmonic in $\{u>0\}$, $u^- \in W^{1,q}(\Omega)$ is $q$-harmonic in $\{u<0\}$ and $\partial\{u>0\}\cap\Omega = \partial\{u<0\}\cap\Omega$. Suppose that $\partial \{u > \e\}\cap\Omega$ and $\partial\{ u< -\delta\} \cap \Omega$ has locally finite perimeter for a.e.\ $\e,\delta>0$. Then 
\begin{equation}\label{eq:fb}
\begin{aligned}
&\lim_{\delta \to 0}(q-1)\int_{\Omega\cap \partial\{u<-\delta\}} |D u|^{q-1} D u\cdot\eta\,dH^{n-1} \\
& = \lim_{\e \to 0} (p-1)\int_{\Omega\cap \partial \{u>\e\}} |D u|^{p-1} D u \cdot \eta\,dH^{n-1},
\end{aligned}
\end{equation}
for any $\eta\in C_0^1(\Omega;\R^n)$.
\end{lem}

\begin{proof}
The first part of the statement follows easily from the fact that $u\in C_{loc}^{0,\sigma}\subset C(\Omega)$ (Theorem \ref{thm:holder}). The second part on the weak formulation \eqref{eq:fb} is also considered as classical. One may prove it by following the computation in \cite[Theorem 2.4]{ACF} with slight modification due to the involvement of the $p$- and $q$-Laplace operator. We skip the detail. 
\end{proof} 

An immediate consequence of the above lemma is that if $\partial\{u>0\}\cap B$ is a $C^1$-graph for some ball $B\subset\Omega$, then 
\begin{equation}\label{eq:fb-re}
(q-1)| D u^-| ^q = (p-1) |D u^+|^p\quad\text{on }\partial\{u>0\}\cap B. 
\end{equation}
Hence, one may Proposition \ref{prop:visc}. 

Next we observe that local minimizers grows at most linearly around vanishing points that can be touched by a ball from either side of the free boundary.

\begin{prop}\label{prop:lin-growth}
Let $u\in W^{1,p\wedge q}(\Omega)$ be a local minimizer , and let $z\in \partial \{u > 0\}\cap \Omega$. Suppose that there is a ball $B\subset \{u >0\} \cup \{ u<0\}$ such that $\partial B\cap \partial\{u>0\} = \{z\}$. Then 
$$
\limsup_{x\to z} \frac{|u(x)|}{|x-z|} < \infty.
$$
\end{prop}

\begin{proof}
Assume without loss of generality that $B \subset \{ u>0\}$. According to Lemma \ref{lem:euler}, $\partial\{u > 0\}\cap\Omega = \partial\{u<0\}\cap\Omega$. Hence, $B$ can be considered as a ball contained in $\{u < 0\}^c$ such that $\partial B\cap\partial \{u<0\} = \{z\}$. Thanks to Theorem \ref{thm:holder} and Lemma \ref{lem:euler}, $u^-$ is a positive $q$-harmonic function in $\{u<0\}$ that vanishes continuously on $\partial\{u<0\}\cap\Omega$. Hence, by Hopf's lemma, 
$$
\limsup\limits_{x\to z} \frac{u^-(x)}{|x-z|} < \infty.
$$ 

Next, we prove that 
$$
\limsup\limits_{x\to z} \frac{u^+(x)}{|x-z|} < \infty.
$$ 
Assume, by way of contradiction, that 
$$
\limsup\limits_{x\to 0} |x|^{-1} u^+(x)= \infty.
$$ 
By Lemma \ref{lem:basic0}, $u^+ \in  L_{loc}^\infty(\Omega)$, so there must exist some small real $r_0>0$ and a sequence $\{r_j \}_{j=1}^\infty\subset(0,r_0)$ of decreasing numbers such that $r_j\to 0$ as $j\to\infty$, and 
\begin{equation}\label{eq:Lip-2}
\sup_{r\in (r_j,r_0)} \left( \frac{1}{r} \sup_{B_r(z)} u^+\right) \leq \frac{1}{r_j} \sup_{B_{r_j}(z)} u^+ = j,
\end{equation} 
for each $j=1,2,\cdots$.

Define an auxiliary function $v_j : B_{r_j^{-1}r_0}\to \R$ by 
$$
v_j(y) = \frac{u^+(r_jy + z)}{jr_j} - \frac{u^-(r_jy + z)}{j^{\frac{p}{q}} r_j}.
$$
Clearly, $v_j\in W^{1,p}(B_{r_j^{-1}r_0}) \cup W^{1,q}( B_{r_j^{-1}r_0}) \cap C_{loc}^{0,\sigma}(B_{r_j^{-1}r_0})$ and, by Lemma \ref{lem:scale}, it is a local minimizer. Moreover, it follows from \eqref{eq:Lip-2} and the assumption $z\in\partial\{u>0\}\cap \Omega$ that
\begin{equation}\label{eq:vj+-sup}
v_j(0) = 0\quad\text{and}\quad \sup_{R\in (1,\frac{r_0}{r_j})} \left( \frac{1}{R} \sup_{B_R} v_j^+\right) \leq \sup_{B_1} v_j^+ = 1.
\end{equation}
Due to the observation above that $\limsup\limits_{x\to z} |x-z|^{-1} u^-(x) < \infty$, we can find some constant $c>0$, independent of $j$, such that 
\begin{equation}\label{eq:vj--sup}
\sup_{R\in (1,\frac{r_0}{r_j})}\left( \frac{1}{R} \sup_{B_R} v_j^-\right) \leq \frac{c}{j^{\frac{p}{q}}}.
\end{equation} 

Combining \eqref{eq:vj+-sup}, \eqref{eq:vj--sup} with Lemma \ref{lem:basic0} and Theorem \ref{thm:holder}, we observe that for each $R\geq 1$, $\{v_j^+\}_{j=j_R}^\infty$ is bounded in $W^{1,p}(B_R)\cap C^{0,\sigma}(B_R)$ and similarly $\{v_j^-\}_{j=j_R}^\infty$ is bounded in $W^{1,q}(B_R)\cap C^{0,\sigma}(B_R)$, where $j_R\geq 1$ is certain large integer. Thus, there is a function $v_0\in W_{loc}^{1,p}(\R^n)\cup W_{loc}^{1,q}(\R^n)\cap C_{loc}^{0,\sigma}(\R^n)$ such that $v_j^+ \to v_0^+$ weakly in $W_{loc}^{1,p}(\R^n)$, $v_j^+\to v_0^-$ weakly in $W_{loc}^{1,q}(\R^n)$ and $v_j \to v_0$ locally uniformly in $\R^n$, along a subsequence. To simplify the notation, we shall continue to denote the convergent subsequence by $\{v_j\}_{j=1}^\infty$. 

In particular, the weak convergence implies that $v_0$ is a local minimizer. In addition, passing to the limit in \eqref{eq:vj+-sup} and utilising the uniform convergence, we observe that 
\begin{equation}\label{eq:v0+}
v_0(0) = 0\quad\text{and}\quad \sup_{R\geq 1} \left( \frac{1}{R} \sup_{B_R} v_0^+\right) \leq \sup_{B_1} v_0^+ = 1.
\end{equation}
Moreover, it follows from \eqref{eq:vj--sup} that $v_0 \geq 0$ in $\R^n$. Thus, $v_0 = v_0^+ \in W_{loc}^{1,p}(\R^n)$, and it implies that $v_0$ minimises $\int_{B_R} |D v|^p \,dy$ over all nonnegative $v\in v_0 + W_0^{1,p}(B_R)$, for each $R\geq 1$ (because $J(w,E) =\int_E |D w|^p\,dy$ for any nonnegative $w\in W^{1,p}(E)$). Therefore, $v_0$ is also a nonnegative $p$-harmonic function in $\R^n$. Thus, the equalities in \eqref{eq:v0+} are incompatible with each other, due to the strong minimum principle. This finishes the proof for $\limsup\limits_{x\to 0} |x|^{-1} u^+(x) < \infty$. 
\end{proof} 

We are ready to prove that local minimizers are viscosity solutions, in the sense of Definition \ref{defn:visc}. 

\begin{proof}[Proof of Proposition \ref{prop:visc}]
Thanks to Theorem \ref{thm:holder}, $u\in C_{loc}^{0,\sigma}(\Omega)\subset C(\Omega)$, for some $\sigma\in(0,1)$ depending only on $n$, $p$ and $q$, and by Lemma \ref{lem:euler}, $u^+\in W^{1,p}(\Omega)$ is a $p$-harmonic function in $\{u > 0\}$ and $u^- \in W^{1,q}(\Omega)$ is a weak harmonic function in $\{u<0\} = \mbox{int}(\{ u\leq 0\})$. Hence, it only remains for us to check the free boundary condition (a) and (b) in Definition \ref{defn:visc}. We shall only present the proof for (a) as the other case follows similarly. 

Let $B\subset \{u>0\}$ be a ball such that $\partial B\cap \partial\{u>0\} = \{z\}$. To simplify the notation, let us consider the case $B = B_1(e_n)$ and $z = 0$ only. Due to the above observation, one may apply \cite[Lemma 4.3]{DK} to both $u^+$ and $u^-$, and then use Proposition \ref{prop:lin-growth} to ensure that there are some (finite) numbers $\alpha > 0$ and $\beta \geq 0$ for which
\begin{equation}\label{eq:visc1}
u^+(x)  \geq \alpha x_n^+ + o(|x_n|)\quad\text{for } x\in B_1(e_n),
\end{equation}
and
\begin{equation}\label{eq:visc2}
u^-(x) \leq \beta x_n^- + o(|x_n|)\quad\text{for }x\in \Omega\setminus B_1(e_n), 
\end{equation} 
as $x\to z$, with equality along every non-tangential domain in each case. Hence, we only need to show that 
\begin{equation}\label{eq:visc3}
(p-1)\alpha^p = (q-1) \beta^q. 
\end{equation}

Consider a linearly scaled version $u_r : r^{-1}  \Omega \to\R$ defined by  
$$
u_r (y) = \frac{u(ry)}{r}.
$$
Thanks to Proposition \ref{prop:lin-growth}, 
$$
\sup_{0<r<r_0}\| u_r \|_{L^\infty(B_{r^{-1} r_0})} \leq c,
$$
for some constant $c>0$. Thus, by Lemma \ref{lem:basic0} and Theorem \ref{thm:holder}, $\{u_r^+\}_{0<r<r_0}$ is bounded in $W^{1,p}(B_R)\cap C^{0,\sigma}(B_R)$ and $\{u_r^-\}_{0<r<r_0}$ is bounded in 
$W^{1,q}(B_R)\cap C^{0,\sigma}(B_R)$, for each $R \in (1,r_0/r_j)$. Therefore, one can find a sequence $r_k\to 0$ of positive numbers and a function $u_0 \in  W_{loc}^{1,p\wedge q}(\R^n) \cap C_{loc}^{0,\sigma}(\R^n)$ such that $u_{r_k}^+ \to u_0^+$ weakly in $W_{loc}^{1,p}(\R^n)$, $u_{r_k}^- \to u_0^-$ weakly in $W_{loc}^{1,q}(\R^n)$ and $u_{r_k} \to u_0$ locally uniformly in $\R^n$. Owing to the weak convergence, $u_0$ is a local minimizer. 

On the other hand, recall that both \eqref{eq:visc1} and \eqref{eq:visc2} hold with equalities in any non-tangential domain. Hence, for each $R\geq 1$ and each $\delta \in (0,1)$, we may find some large integer $k_0\geq 1$ such that 
$$
\sup_{y\in B_R\cap K_\delta} |u_{r_k}^+ (y) - \alpha y_n| + \sup_{y\in B_R \cap (-K_\delta)} |u_{r_k}^- (y) - \beta y_n| \leq \frac{1}{k}, 
$$
for all $k\geq k_0$, where $K_\delta = \{y\in\R^n: y_n \geq  \delta |y|\}$. Thus, sending $k\to\infty$, utilising the uniform convergence of $u_{r_k} \to u_0$ in $B_R$, and then letting $R\to \infty$ and $\delta\to 0$, we obtain that 
$$
u_0^+ (y) = \alpha y_n^+ - \beta y_n^-\quad\text{for all }y\in\R^n. 
$$
In particular, $D u_0 = \alpha e_n$ on $\{y_n > 0\}$, $D u_0 = -\beta e_n$ on $\{y_n < 0\}$, $\partial\{u_0> \alpha \e\} = \{y_n = \e\}$, and $\partial\{ u_0 < -\beta\e\} = \{y_n = -\e\}$. As $u_0$ being a local minimizer whose level sets being hyperplanes, we may invoke Lemma \ref{lem:euler} to derive that 
$$
\begin{aligned}
\beta^p \int_{B_R\cap\{ y_n =0\}} e_n\cdot \eta \,dH^{n-1} &= \lim_{\e\to 0}\int_{B_R\cap \{y_n = -\e\}} |D u_0|^{q-1}D u_0 \cdot \eta \,dH^{n-1} \\
&= \frac{p-1}{q-1}\lim_{\e\to 0}  \int_{B_R\cap \{y_n = \e\}} |D u_0|^{p-1} D u\cdot \eta \,dH^{n-1}\\
 &= \frac{p-1}{q-1} \alpha^p \int_{B_R\cap \{y_n = 0\}} e_n\cdot \eta\,dH^{n-1},
\end{aligned}
$$
for any $\eta\in C_0^1(B_R;\R^n)$. This proves the desired relation \eqref{eq:visc3} between $\alpha$ and $\beta$.
\end{proof}

Now let us turn our attention  to  Proposition \ref{prop:visc}. The proof of this proposition will repeat (almost verbatim)  works  of J. Lewis and K. Nystr\"om \cite{LN1,LN2}
 concerning the two-phase Bernoulli problem with single $p$-Laplace operator. 
These authors also remark  that it is straightforward to extend  their argument to the case of the 
$p$-Laplace operator on one side and the $q$-Laplace operator on the other side; see \cite[Page 108]{LN2}.

One may check the argument in \cite{LN1,LN2} in detail, and verify that their remark is true. Nevertheless, we shall not do it here, as it would only reproduce lengthy argument, of an already known techniques and approach.
 Instead, we shall outline the argument briefly, and point out the differences that need to be addressed, leave some detail in Appendix \ref{appendix:flat}, and conclude the section. The writing here (including the appendix) will also be similar, to some extent, with \cite{Fel1}, which extends to the anisotropic case the seminal papers \cite{C1,C2} on the viscosity method for single standard Laplace operator.

\begin{proof}[Proof of  Proposition \ref{prop:flat-fb} ]  
Let $u\in C(\Omega)$ be a viscosity solution, in the sense of Definition \ref{defn:visc}. There are two points that we need to address before we repeat the proof of \cite[Theorem 1]{LN1} and \cite[Theorem 1]{LN2}. Firstly, the definition of viscosity solutions in \cite{LN1,LN2} (see \cite[Definition 1.4]{LN1}) is different from Definition \ref{defn:visc} in the sense that here we impose two different operators $\Delta_p$, $\Delta_q$, on different phases. This affects \cite[Lemma 3.5, Lemma 3.22]{LN1} and \cite[Lemma 3.7]{LN2}, which are taken care of by Lemma \ref{lem:flat}, Lemma \ref{lem:fam1} and respectively Lemma \ref{lem:fam2}. Secondly, in the hypothesis of \cite[Theorem 1]{LN1} and \cite[Theorem 1]{LN2}, the function $G:[0,\infty)\to [0,\infty)$ for the balance equation, $\alpha = G(\beta)$, between the slopes $\alpha,\beta$ of linear asymptotic developments should satisfy that $G > 0$ is strictly increasing in $[0,\infty)$ and that $s\mapsto s^{-N}G(s)$ is decreasing in $[0,\infty)$ for some $N>0$. In our definition (see Definition \ref{defn:visc}), $G(s) = (\frac{p-1}{q-1} s^q)^{1/p}$ for any $s\geq 0$, so $G\geq 0$ (with equality at $s=0$) is strictly increasing and $s\mapsto s^{-N}G(s)$ is decreasing with any $N > \frac{q}{p}$. The fact that $G(0) = 0$ in our definition does not affect the argument of \cite[Theorem 1]{LN1} and \cite[Theorem 1]{LN2}, since we always have $\alpha>0$ and $\beta>0$, while the aforementioned references also cover the case where $\alpha > 0$ but $\beta \geq 0$.\footnote{As a matter of fact, this makes the argument here to be less complicated, since we do not need to consider the case when $\alpha > 0$ but $\beta = 0$ (i.e., $u^+$ is non-degenerate but $u^-$ is degenerate), which is a challenging part in the analysis on two-phase Bernoulli problems.} 
We  leave  the details  to the reader, and finish the proof here. 
\end{proof} 

Let us finish the section with the proof for our main result. 

\bigskip
\noindent
{\bf Proof of Theorem \ref{thm:flat}:}
Due to Proposition \ref{prop:visc}, a local minimizer is a viscosity solution. Hence, Proposition \ref{prop:flat-fb} applies to any flat free boundary point of the local minimizer. This finishes the proof.


\section{Structure of  the Free Boundary (Proof of Theorem \ref{thm:nonflat})}\label{section:nonflat}

Here we shall study the property of the measure $\Delta_p u^+$, and how small the non-flat part of the free boundary of a local minimizer $u$ can be with respect to this measure. Note that as $u^+ \in W^{1,p}(\Omega)$ is a $p$-harmonic function in $\{u>0\}$ continuously vanishes to the free boundary $\partial\{u>0\}\cap\Omega$, we can define $\Delta_p u^+$ as a nonnegative Radon measure, by setting
$$
\int_\Omega \phi \Delta_p u^+ = -\int_\Omega |D u^+|^{p-2} D u^+ \cdot D \phi\,dx,
$$ 
for any $\phi\in C_c^\infty(\Omega)$. Clearly, $\Delta_p u^+$ is supported on $\partial\{u>0\}\cap\Omega$. Our main assertion of this section is Theorem \ref{thm:nonflat}. 

This type of result was considered by Andersson and Mikayelyan \cite{AM}, concerning the free boundaries of weak solutions to an elliptic anisotropic problem. However, new challenges arise from the presence of the $p$-Laplacian on the one side and the $q$-Laplacian on the other. The situation becomes very different when we attempt to obtain a nontrivial blowup limit only with a doubling condition from one side (see Lemma \ref{lem:blowup}). The main difficulty lies in the fact that the correct scaling, by which the scaled version continues to be a local minimizer, requires comparability between $\norm{u^-}_{L^\infty(B_r(z))}$ and $r^{1-p/q}\norm{u^+}_{L^\infty(B_r(z))}^{p/q}$, rather than $\norm{u^+}_{L^\infty(B_r(z))}$. Let us remark that such an issue does not appear for standard anisotropic problems, where two different operators are of the same order. 

Let us begin with some basic observations regarding the measure $\Delta_p u^+$ for a weak $p$-subsolution $u^+$. The first is the natural growth estimate. As a byproduct, we observe that $\Delta_p u^+$ is a locally finite measure. 

\begin{lem}\label{lem:msr-growth}
Let $u^+\in W^{1,p}(\Omega)$ be a weak $p$-subsolution. Then for any $z\in\Omega$ and any $r\in(0,\frac{1}{4}\dist(z,\partial\Omega))$, 
\begin{equation}\label{eq:msr-growth}
\frac{\Delta_p u^+(B_r(z))}{r^{n-p}} \leq C \norm{u^+}_{L^\infty(B_{2r}(z))}^{p-1},
\end{equation} 
where $C>0$ depends only on $n$, $p$ and $q$. 
\end{lem} 

\begin{proof}
For the sake of simplicity, let us consider $z = 0$ and $\Omega = B_1$. This inequality follows easily from the Caccioppoli inequality and the local maximum principle of weak $p$-subsolutions. Let $r\in(0,\frac{1}{4})$ be given, so that $\overline{B_{2r}}\subset B_{1/2}$. Let $\phi\in C_c^\infty(B_{2r})$ be a smooth cutoff function such that $\phi = 1$ on $\overline{B_r}$, $0\leq \phi\leq 1$ in $B_{2r}$ and $\supp \phi \subset B_{3r/2}$. Then as $\Delta_p u^+$ being a positive Radon measure, 
$$
\begin{aligned}
\Delta_p u^+(B_{2r}) & \leq \int_{B_{3r/2}} \phi \,\Delta_p u^+  \\
& \leq \int_{B_{3r/2}} |D u^+|^{p-1} |D \phi|\,dx   \leq \frac{c_1}{r} \int_{B_{3r/2}} |D u^+|^{p-1}\,dx \\
& \leq  c_2 r^{n-1} \left( \frac{1}{r^n}\int_{B_{3r/2}} |D u^+|^p\,dx \right)^{1- 1/p}  \leq  c_3 r^{n-p} \norm{u^+}_{L^\infty(B_{2r})}^{p-1},
\end{aligned}
$$  
where in the derivation of the last two lines we applied Jensen's inequality and the Caccioppoli inequality. Here $c_1$, $c_2$ and $c_3$ are constants depending at most on $n$ and $p$. 
\end{proof} 

The next lemma is an analogue of Lemma 4.3 in \cite{AM}, which asserts that the set of highly  degenerate points are small with respect to measure $\Delta_p u^+$. 

\begin{lem}\label{lem:msr-deg} 
Let $u^+ \in W^{1,p}(\Omega)$ be a weak $p$-subsolution, and consider 
$$
\Gamma_k :=\left \{ z\in \Omega: \limsup\limits_{r\to 0} \frac{\Delta_p u^+(B_r(z))}{r^k} < \infty \right\}. 
$$
Then $\Delta_p u^+(\Gamma_k) = 0$ for any $k>n$. 
\end{lem} 

\begin{proof}
Let $M\geq 1$ and $\e\in(0,1)$ be given, and set 
$$
\Gamma_{k,M,\e} := \left\{z\in\Omega_\e: \Delta_p u^+(B_r(z)) \leq M r^k,\,\forall r\in(0,\e) \right\},
$$
where $\Omega_\e := \{ z\in\Omega:  \dist(z,\partial\Omega) \geq 2\e\}$. We shall prove that $\Delta_p u^+(\Gamma_{k,M,\e}) = 0$, for any $k>n$. This suffices to justify the assertion of the lemma, since $\Gamma_k = \bigcup_{\ell,m=1}^\infty \Gamma_{k,\ell,\frac{1}{m}}$. 

Note that if $z_i\to z \in B_{1-2r}$, then 
$$
\lim\limits_{i\to\infty} \int_{B_r(z_i)} |D u^+|^{p-2} D u^+\cdot D \phi\,dx = \int_{B_r(z)} |D u^+|^{p-2} D u^+\cdot D \phi\,dx,
$$
for any nonnegative $\phi\in C_0^\infty(\Omega)$. Hence, $\Gamma_{k,M,\e}$ is a closed set, thus a compact set. Therefore, for each $r\in(0,\e)$, we can choose $\{z_i\}_{i=1}^N\subset \Gamma_{k,M,\e}$, for some constant $N \leq c_0r^{-n}$ with $c_0>0$ depending only on $n$, such that $\Gamma_{k,M,\e}\subset \bigcup_{i=1}^N B_r(z_i)$. According to  Vitali's covering lemma, we can also choose $\{z_i\}_{i=1}^N$ in such a way that $B_{r/5}(z_i)\cap B_{r/5}(z_j) = \emptyset$ for any distinct pair $(i,j)$ of indices from $\{1,\cdots,N\}$, and each $\#\{j\in\{1,\cdots,N\}: B_r(z_i)\cap B_r(z_j) \neq \emptyset\} \leq c_0$, where $c_0$ is a constant depending only on $n$. Then from the fact that $z_i \in \Gamma_{k,M,\e}$, and that $\bigcup_{i=1}^N B_r(z_i)\subset \Omega$, we deduce that for any $r\in(0,\e)$, 
$$
\Delta_p u^+(\Gamma_{k,M,\e}) \leq \sum_{i=1}^N \Delta_p u^+(B_r(z_i)) \leq MNr^k \leq \frac{c_0 r^{k-n} M|\Omega|}{|B_1|},
$$ 
so letting $r\to 0$ yields that $\Delta_p u^+(\Gamma_{k,M,\e}) = 0$, as desired. 
\end{proof} 

From the lemma above, we deduce that $u$ satisfies the doubling property almost everywhere with respect to the measure $\Delta_p u^+$, at least along some subsequences. 

\begin{lem}\label{lem:msr-double}
Let $u^+ \in W^{1,p}(\Omega)\cap C(\Omega)$ be a weak $p$-subsolution, and define 
$$
F := \left\{ z\in \partial\{u>0\}\cap \Omega: \liminf\limits_{r\to 0} \frac{\norm{u^+}_{L^\infty(B_{2r}(z))}}{\norm{u^+}_{L^\infty(B_r(z))}} = \infty \right\}. 
$$
Then $\Delta_p u^+(F) = 0$. 
\end{lem} 

\begin{proof}
Let $\Gamma_k$ be as in Lemma \ref{lem:msr-deg}, so that $\Delta_p u^+(\Gamma_k) = 0$ for any $k>n$, hence for $k = n +1$. By the definition of $\Gamma_{n+1}$, we have, for any $z\in\partial\{u>0\}\cap \setminus \Gamma_{n+1}$, that there exists a sequence $\{\rho_k\}_{k=1}^\infty$ of positive real numbers such that $\rho_k\to 0$ as $k\to\infty$, and 
$$
k \leq \frac{ \Delta_p u^+(B_{\rho_k}(z))}{\rho_k^{n+1}} \leq \frac{C \| u^+ \|_{L^p(B_{2\rho_k}(z))}^{p-1}}{\rho_k^{p+1}},
$$
for each $k=1,2,\cdots$, where the second inequality follows from \eqref{eq:msr-growth}. Thus, we obtain   
$$
\limsup_{r\to 0} \left( r^{- \frac{p+1}{p-1}} \norm{u^+}_{L^\infty(B_r(z))} \right) = \infty.
$$ 
From the above observation, we can choose a sequence $\{r_j\}_{j=1}^\infty\subset(0,r_0)$, with $r_0< \frac{1}{2} \dist(z,\partial \Omega)$, of positive real numbers, decreasing to zero, such that
$$
\sup_{r\in (r_j,2r_0)} \left( r^{-\frac{p+1}{p-1}} \norm{u^+}_{L^\infty(B_r(z))} \right) \leq r_j^{-\frac{p+1}{p-1}} \norm{u^+}_{L^\infty(B_{r_j}(z))}, 
$$ 
from which it follows immediately that $z\in \partial\{u>0\}\setminus F$. Therefore, we have proved that $F\subset \Gamma_{n+1} $, showing that $F$ has null $\Delta_p u^+$-measure. 
\end{proof} 

Note that the last three lemmas are true for any weak $p$-subsolution. Henceforth, we shall study stronger properties for local minimizers. First, let us observe that the sequence $\{u_{z,r}\}_{r>0}$ of $u$, with $u_{z,r}$ defined by
\begin{equation}\label{eq:uzr}
u_{z,r}(x) = \frac{u^+(2rx+z)}{\norm{u^+}_{L^\infty(B_r(z))}} - \frac{u^-(2rx+z)}{(2r)^{1-\frac{p}{q}} \norm{u^+}_{L^\infty(B_r(z))}^{\frac{p}{q}}},
\end{equation}
at a free boundary point $z$, where $u^+$ satisfies a doubling condition along a sequence, contains a subsequence that converges to a non-trivial local minimizer. 

\begin{lem}\label{lem:blowup}
Let $u\in W^{1,p\wedge q}(\Omega)$ be a local minimizer, and let $z\in\partial\{u>0\}\cap\Omega$ and $r \in (0,\frac{1}{2}\dist(z,\partial\Omega))$ be arbitrary. If 
$$
\liminf\limits_{r\to 0} \frac{\norm{u}_{L^\infty(B_{2r}(z))}}{\norm{u}_{L^\infty(B_r(z))}} < \infty,
$$
then there exists a nontrivial local minimizer $u_{z,0} \in W_{loc}^{1,p\wedge q}(B_1)$ such that $u_{z,r} \to u_{z,0}$ locally uniformly on $B_1$, $D u_{z,r}^+ \to D u_{z,0}^+$ weakly in $L_{loc}^p(B_1)$ and $D u_{z,r}^- \to D u_{z,0}^-$ weakly in $L_{loc}^q(B_1)$, along a subsequence as $r\to 0$. 
\end{lem} 

\begin{proof}
By the assumption, we can find a sequence $\{r_j\}_{j=1}^\infty\subset(0,\frac{1}{2})$ of positive real numbers decreasing to $0$ and a constant $\beta \in (0,1)$ such that  
\begin{equation}\label{eq:double}
\norm{u^+}_{L^\infty(B_{2r_j}(z))} \leq \frac{1}{\beta} \norm{u^+}_{L^\infty(B_{r_j}(z))}. 
\end{equation} 
Define $u_j = u_{z,r_j} \in W^{1,p\wedge q}(B_1)\cap C(\overline{B_1})$, where $u_{z,r}$ is given as in \eqref{eq:uzr}. By Lemma \ref{lem:scale} and \eqref{eq:double}, $u_j$ is a local minimizer such that 
\begin{equation}\label{eq:uj+-linf}
\norm{u_j^+}_{L^\infty(B_{1/2})} = 1,\quad \norm{u_j^+}_{L^\infty(B_1)} \leq \frac{1}{\beta}. 
\end{equation}

For the rest of the proof, we shall denote by $c$ a generic constant depending at most on $n$, $p$, $q$ and $\beta$; in particular, $c$ may vary at each occurrence. 

Thanks to \eqref{eq:uj+-linf}, it follows from Lemma \ref{lem:basic0}, Caccioppoli inequality and Theorem \ref{thm:holder} that $\{ u_j^+ \}_{j=1}^\infty$ is a bounded sequence in $W_{loc}^{1,p}(B_1)\cap C_{loc}^{0,\sigma}(B_1)$, where $\sigma \in(0,1)$ depends only on $n$, $p$ and $q$; it is noteworthy that the boundedness can be deduced without any estimate on $\norm{u_j^-}_{L^q(B_1)}$. In particular, for each $r_0 \in (\frac{1}{2},1)$, we have 
\begin{equation}\label{eq:uj+-Ca}
\| D u_j^+ \|_{L^p(B_r)}  \leq \frac{c}{(1-r_0)^{\frac{n}{p}}}, \quad [ u_j^+ ]_{C^{0,\sigma}(B_r)} \leq \frac{c}{(1-r_0)^{\frac{n}{p} + \sigma}}. 
\end{equation}

Due to the fact that $\norm{u_j^+}_{L^\infty(B_{1/2})} = 1$ and $u_j^+\in C(B_1)$, there exists a point $x_j\in \overline{B_{1/2}}$ such that $u_j^+(x_j) = 1$. By \eqref{eq:uj+-Ca}, there is some constant $\rho\in(0,\frac{1}{8})$ such that $\inf\limits_{B_\rho(x_j)} u_j^+ \geq 0$; here $\rho$ depends only on $c$, $n$, $p$ and $\sigma$, hence independent on $j$. This implies that $\{ u_j < 0\}\subset B_r \setminus B_\rho(x_j)$ for any $r\in(\frac{3}{4},1)$, which in turn yields that 
\begin{equation}\label{eq:uj+-msr}
\frac{| \{ u_j < 0 \} \cap B_{r_0} |}{|B_{r_0}|} \leq 1 - \frac{|B_\rho (x_j)|}{|B_{r_0}|} \leq \gamma, 
\end{equation} 
for some constant $\gamma_n\in(0,1)$; note that since $\rho\in(0,\frac{1}{8})$, $x_j\in \overline{B_1}$ and $r_0\in(\frac{3}{4},1)$, we have $B_\rho(x_j) \subset B_{5/8}$, so $\gamma$ depends only on $n$. Hence, the hypothesis of Lemma \ref{lem:holder3} is verified, from which we obtain that for each $r_0\in(\frac{3}{4},1)$, 
\begin{equation}\label{eq:uj--Linf}
\norm{u_j^-}_{L^\infty(B_{r_0})} \leq \frac{c}{(1-r_0)^\mu} \norm{u_j^+}_{L^\infty(B_1)}^{\frac{p}{q}} \leq \frac{c}{(1-r_0)^\mu}, 
\end{equation} 
where $\mu>0$ depends only on $n$, $p$ and $q$. As $r_0\in(\frac{3}{4},1)$ being arbitrary, we may follow the derivation of \eqref{eq:uj+-Ca} and deduce that $\{ u_j^-\}_{j=1}^\infty$ is a bounded sequence in $W_{loc}^{1,q}(B_1)\cap C_{loc}^{0,\sigma}(B_1)$. 

Owing to the boundedness of $\{u_j^+\}_{j=1}^\infty$ and $\{u_j^-\}_{j=1}^\infty$ in $W_{loc}^{1,p}(B_1)\cap C_{loc}^{0,\sigma}$ and respectively $W_{loc}^{1,q}(B_1)\cap C_{loc}^{0,\sigma}(B_1)$, we can extract a subsequence $\{u_{j_k}\}_{k=1}^\infty$ of $\{u_j\}_{j=1}^\infty$ and find a function $u_0 \in W_{loc}^{1,p\wedge q}(B_1)\cap C_{loc}^{0,\sigma}(B_1)$ such that $u_{j_k} \to u_0$ locally uniformly on $B_1$, $D u_{j_k}^+ \to D u_0^+$ weakly in $L_{loc}^p(B_1)$ and $D u_{j_k}^- \to D u_0^-$ weakly in $L_{loc}^q(B_1)$, as $k\to \infty$. Due to the weak convergence, $u_0$ is a local minimizer. Moreover, utilising the uniform convergence, we may deduce from \eqref{eq:uj+-linf} that 
\begin{equation}\label{eq:u0}
\norm{u_0^+}_{L^\infty(B_{1/2})} = 1.
\end{equation}
This proves that $u_0$ is a non-trivial local minimizer, as desired. 
\end{proof}

We are in a position to assert that the free boundary of a local minimizer $u$ is flat a.e.\ with respect to the measure $\Delta_p u^+$. Analogously as with \cite{DK}, for each point $z\in\partial\{u>0\}\cap \Omega$, each direction $\nu\in \partial B_1$ and each scale $r\in(0,\dist(z,\partial \Omega))$, let $h(z,\nu,r)$ denote the minimum height of $\partial\{u>0\}\cap B_r(z)$, i.e.,
\begin{equation}\label{eq:h}
h(z,\nu,r) = \sup\{ |\inn{x-z}{\nu} | : x\in \partial\{u>0\}\cap B_r(z)\},
\end{equation}
and then define 
\begin{equation}\label{eq:hb}
\bar h (z,r) = \inf_{\nu\in \partial B_1} h(z,\nu,r). 
\end{equation}

\begin{prop}\label{prop:flat} 
Let $u \in W^{1,p\wedge q}(\Omega)$ be a local minimizer, and define 
$$
E := \left\{z\in\partial \{u>0\}\cap \Omega: \liminf\limits_{r\to 0} \frac{\bar h(z,r)}{r} > 0 \right\}.
$$
Then $\Delta_p u^+(E) = 0$. 
\end{prop} 

\begin{proof}
Let $F$ be the set as in Lemma \ref{lem:msr-double}, and for each $k,\ell=1,2,\cdots$, define 
$$
E_{k,\ell} := \left\{z\in E:  \bar h(z,r) \geq \frac{r}{k},\,\forall r\in(0,\frac{1}{2\ell}) \right\}\cap \left\{z\in\Omega: \dist(z,\partial E) >\frac{1}{\ell} \right\}. 
$$
We claim that 
$$
\Delta_p u^+( E_{k,\ell}\setminus F) = 0,\quad\forall k,\ell=1,2,\cdots, 
$$
which along with $\Delta_p u^+ (F) = 0$ proves that $\Delta_p u^+( E_{k,\ell})= 0$. Consequently, we will be able to conclude from $E = \bigcup_{k,\ell=1}^\infty E_{k,\ell}$ that $\Delta_p u^+( E) = 0$. 

Fix a pair $(k,\ell)$ of positive integers. To prove that $\Delta_p u^+(E_{k,\ell}\setminus F) = 0$, it suffices to verify that for each $z\in E_{k,\ell}\setminus F$, there exists some $\delta \in (0,1)$ such that 
\begin{equation}\label{eq:claim-flat}
\liminf_{r\to 0} \frac{ \Delta_p u^+(E_{k,\ell} \cap B_r(z)) }{ \Delta_p u^+(B_r(z)) } \leq 1-\delta.
\end{equation} 
This implies that the density points of $E_{k,\ell}\setminus F$ constitute a set of null $\Delta_p u^+$-measure, so we must have $\Delta_p u^+( E_{k,\ell}\setminus F) = 0$.  

Fix a point $z\in E_{k,\ell}\setminus F$. Let $u_{z,r} \in W^{1,p\wedge q}(B_1) \cap L^\infty(B_1)$ be defined as in \eqref{eq:uzr}. By Lemma \ref{lem:blowup}, we may find a sequence $\{r_j\}_{j=1}^\infty$ of positive real numbers decreasing to $0$ such that $u_{z,r_j} \to u_{z,0}$ locally uniformly on $B_1$, $D u_{z,r_j}^+ \to D u_{z,0}^+$ weakly in $L_{loc}^p(B_1)$ and $D u_{z,r_j}^- \to D u_{z,0}^-$ weakly in $L_{loc}^q(B_1)$, for certain nontrivial local minimizer $u_{z,0}\in W_{loc}^{1,p\wedge q}(B_1)\cap C_{loc}^{0,\sigma}(B_1)$. Note also from the assumption $z\in E_{k,\ell}\setminus F\subset \partial\{u>0\}\cap \Omega$ that $u_{z,r_j}(0) = 0$ for all $j=1,2,\cdots$, so the uniform convergence implies that $u_{z,0}(0) = 0$. 

According to Lemma \ref{lem:euler}, $u_{z,0}^+ \in W_{loc}^{1,p}(B_1)$ is $p$-harmonic in $\{u_{z,0}>0\}$ and $u_{z,0}^-\in W_{loc}^{1,q}(B_1)$ is $q$-harmonic in $\{u_{z,0}<0\}$.
 Thus, it follows from the strong maximum principle and $u_{z,0}(0) = 0$ that $u_{z,0}^+\not \equiv  0$ and $u_{z,0}^- \not \equiv  0$ in $B_r$ for any $r\in(0,1)$. 

Since $u_{z,0}^+ \not \equiv 0$ in $B_{1/2}$ 
and $u_{z,0} \in C_{loc}^{0,\sigma}(B_1)\subset C(B_1)$, we can find a ball $B\subset \{u_{z,0}>0\}\cap B_{1/2}$ such that $\partial B\cap \partial \{u_{z,0}>0\}\cap B_{1/2} \neq \emptyset$. Moreover, we may assume without loss of generality that $\{x_0\} = \partial B\cap \partial \{u_{z,0}>0\}\cap B_{1/2}$. According to Proposition \ref{prop:visc}, $u_{z,0}$ is a viscosity solution in the senes of Definition \ref{defn:visc}, so it follows from Condition (a) in Definition \ref{defn:visc} that 
\begin{equation}\label{eq:u0-twoplane}
\frac{u_{z,0}(\rho y + x_0)}{\rho} \to \alpha \inn{y}{\nu}^+ - \beta \inn{y}{\nu}^-,\quad\text{as }\rho\to 0,
\end{equation} 
uniformly for all $y\in \overline{B_1}$, where $\nu\in\partial B_1$ is the inward unit normal to $\partial B_1$ at $x_0$, and $(\alpha,\beta)$ is a pair of positive constants satisfying $(p-1)\alpha^p = (q-1)\beta^q$.

Due to the uniform convergence in \eqref{eq:u0-twoplane}, one can find a sufficiently small $\rho_0\in(0,\frac{1}{16\ell})$ such that $\partial\{u_{z,0}>0\}\cap B_{8\rho_0}(x_0) \subset \{x\in \R^n: |\inn{x-x_0}{\nu}| \leq \frac{\rho_0}{16k}\}$. Then for any $\xi\in\partial\{u_{z,0}>0\}\cap B_{2\rho_0}(x_0)$, we have 
\begin{equation}\label{eq:flat-u0}
\partial\{u_{z,0}>0\}\cap B_{2\rho_0}(\xi)\subset\left\{x\in\R^n: |\inn{x-\xi}{\nu}| \leq \frac{\rho_0}{8k}\right\}.
\end{equation} 

Recall from $x_0 \in B_{1/2}$ and $\rho_0< \frac{1}{16\ell} \leq \frac{1}{16}$ that $B_{2\rho_0}(x_0)\subset B_{5/8}$. Also recall that $u_{z,r_j} \to u_{z,0}$ locally uniformly on $B_1$, which along with \eqref{eq:flat-u0} implies the following: for any large integer $j\geq 1$, there exist $x_j\in\partial \{u_{z,r_j}>0\}\cap B_{9/16}$ and $\nu_j\in\partial B_1$ such that for any $\xi\in\partial \{u_{z,r_j}>0\}\cap B_{\rho_0}(x_j)$, we have 
\begin{equation}\label{eq:flat-uj}
\partial\{u_{z,r_j}>0\}\cap B_{\rho_0}(\xi)\subset\left\{x\in\R^n: |\inn{x-\xi}{\nu_j}| \leq \frac{\rho_0}{4k}\right\}.
\end{equation} 
In view of \eqref{eq:uzr}, we may rewrite \eqref{eq:flat-uj} in terms of $u$ as that for any $\zeta\in \partial\{u>0\}\cap B_{2\rho_0r_j} (z + 2r_j x_j)$, 
$$
\partial\{u>0\}\cap B_{2\rho_0 r_j}(\zeta)\subset\left\{x\in\R^n: |\inn{x-\zeta}{\nu_j}| \leq \frac{\rho_0r_j}{2k}\right\}.
$$
This implies that 
$$
\bar h(\zeta,r_j) \leq \bar h(\zeta,2r_j) \leq \frac{\rho_0r_j}{2k},\quad\text{for all }\zeta \in \partial \{u>0\}\cap B_{2\rho_0r_j}(z + 2r_jx_j),
$$
where $\bar h(\zeta,r)$ is defined as in \eqref{eq:hb}. Therefore, 
\begin{equation}\label{eq:Ekl}
E_{k,\ell} \cap B_{r_j}(z) = E_{k,\ell} \cap (B_{r_j}(z) \setminus B_{2\rho_0r_j}(z + 2r_jx_j)).
\end{equation}

Henceforth, we shall focus on finding some $\delta\in(0,1)$ such that 
\begin{equation}\label{eq:claim-flat-2}
\liminf_{j\to \infty} \frac{\Delta_p u^+(B_{2\rho_0 r_j}(z + 2r_jx_j)) }{\Delta_p u^+(B_{r_j}(z))} \geq \delta.
\end{equation} 
Note that if \eqref{eq:claim-flat-2} is true, then our initial claim \eqref{eq:claim-flat} follows immediately from \eqref{eq:Ekl}.

Thanks to the uniform convergence in \eqref{eq:u0-twoplane}, $\{\rho^{-1}u_{z,0}(\rho \cdot + x_0)\}_{0<\rho<1}$ is a bounded family in $L^\infty(B_1)$. As $u_{z,0}$ being a local minimizer, it follows from Lemma \ref{lem:basic0} (that $u_{z,0}^+$ is a weak $p$-subsolution) and the Caccioppoli inequality that $\{\rho^{-1} u_{z,0}^+(\rho\cdot + x_0)\}_{0<\rho<1}$ is a bounded family in $W_{loc}^{1,p}(B_1)$. Hence, for any subsequence $\{w_j\}_{j=1}^\infty$ of $\{\rho^{-1} u_{z,0}^+(\rho\cdot + x_0)\}_{0<\rho<1}$, there exists some $w_0 \in W_{loc}^{1,p}(B_1)$ such that $w_j\to w_0$ weakly in $W_{loc}^{1,p}(B_1)$. However, the uniform convergence in \eqref{eq:u0-twoplane} implies that $w_0 = \alpha \inn{\cdot}{\nu}^+$ in $B_1$, regardless of the choice of the subsequence $\{w_j\}_{j=1}^\infty$. Thus, $D u_{z,0}^+(\rho \cdot + x_0) \to \alpha \nu H(\inn{\cdot}{\nu})$ weakly in $L_{loc}^p(B_1)$, where $H$ is the Heaviside function on the real line. 

Hence, denoting by $\Pi$ the hyperplane $\{y\in\R^n: \inn{y}{\nu} = 0\}$, and choosing a cutoff function $\phi\in C_c^\infty(B_1)$ with $\supp \phi \subset B_{3/4}$, $0\leq \phi\leq 1$ in $B_1$ and $\phi = 1$ on $\overline{B_{1/2}}$, we have 
$$
\begin{aligned} 
& \lim_{\rho \to 0} \frac{1}{\rho^{n-1}} \int_{B_\rho(x_0)} |D u_{z,0}^+(x)|^{p-2} D u_{z,0}^+ (x) \cdot D \phi\left(\frac{x - x_0}{\rho}\right)dx  \\
& = \alpha \int_{B_1} \nu\cdot D \phi(y) H(\inn{y}{\nu})\,dy = -\alpha \int_{B_1 \cap \Pi} \phi\,dH^{n-1}  \leq - \alpha \lambda, 
\end{aligned} 
$$ 
with $\lambda = H^{n-1} (B_{1/2}\cap \Pi) > 0$; here the negative sign follows from the fact that $\nu$ is the inward unit normal to $\Pi$. Therefore, taking the above $\rho_0$ smaller if necessary, we obtain that 
\begin{equation}\label{eq:msr-u0}
\Delta_p u_{z,0}^+ ( B_{\rho_0/2}(x_0) ) \geq \frac{\alpha\lambda\rho_0^{n-1}}{2^n}.  
\end{equation} 

Let us also recall that $x_j\in \partial\{u_{z,r_j}>0\}\cap B_{5/8}$, which appears in \eqref{eq:claim-flat-2}, was chosen such that $x_j \to x_0$ as $j\to \infty$. Now recalling that $D u_{z,r_j}^+ \to D u_{z,0}^+$ weakly in $L_{loc}^p(B_1)$, and $B_{3\rho_0/4}(x_0) \subset B_{\rho_0}(x_j) \subset B_1$ for all sufficiently large $j$, we see that
\begin{equation}\label{eq:msr-uj}
\Delta_p u_{z,r_j}^+(B_{\rho_0}(x_j)) \geq \frac{1}{2} \Delta_p u_{z,0}^+ ( B_{\rho_0/2}(x_0)) \geq \frac{\alpha\lambda \rho_0^{n-1}}{2^{n+2}}.
\end{equation} 
Utilising the growth estimate (Lemma \ref{lem:msr-growth}) of the measure $\Delta_p u_{z,r_j}^+$, and the fact that $u_{z,r_j}\to u_{z,0}$ locally uniformly in $B_1$, we observe that 
\begin{equation}\label{eq:msr-uj-2}
\Delta_p u_{z,r_j}^+(B_{1/2}) \leq C\norm{u_{z,r_j}^+}_{L^\infty(B_{3/4})} \leq C\Lambda,
\end{equation} 
for all sufficiently large $j$, where $\Lambda = 2\norm{u_{z,0}^+}_{L^\infty(B_{3/4})}$. Combining \eqref{eq:msr-uj} together with \eqref{eq:msr-uj-2}, we observe that
\begin{equation}
\frac{\Delta_p u^+(B_{2\rho_0 r_j}(z + 2r_jx_j)) }{\Delta_p u^+(B_{r_j}(z))} = \frac{\Delta_p u_{z,r_j}^+(B_{\rho_0}(x_j))}{\Delta_p u_{z,r_j}^+(B_{1/2})} \geq \frac{\alpha\lambda \rho_0^{n-1}}{2^{n+2}C\Lambda},
\end{equation} 
for any large $j$, which proves \eqref{eq:claim-flat-2} with $\delta = \frac{\alpha\lambda \rho_0^{n-1}}{2^{n+2}C\Lambda}>0$, as desired. 
This completes the proof. 
\end{proof} 

Let us finish the section with the proof of Theorem \ref{thm:nonflat}.

\bigskip
\noindent
{\bf Proof of Theorem \ref{thm:nonflat}}
Let $E$ be the set as in Proposition \ref{prop:flat}, and choose any $z\in (\partial\{u>0\}\cap\Omega)\setminus E$. Since $\liminf\limits_{r\to 0} \frac{\bar h(z,r)}{r} = 0$, we can find a sufficiently small $r_0\in(0,\frac{1}{2}\dist(z,\partial\Omega))$ such that $\bar h(z,r_0) \leq \e r_0$, where $\e\in(0,1)$ is the small (flatness) constant chosen as in Theorem \ref{thm:flat}. Then with $\sigma$ and $\rho$ as in Theorem \ref{thm:flat}, we obtain that $\partial\{u>0\}\cap B_{\rho r_0}(z)$ is a $C^{1,\sigma}$-graph, since $u$ as a local minimizer it  is also a viscosity solution, according to Proposition \ref{prop:visc}. Since $\Delta_p u^+(E) = 0$, the first part of Theorem \ref{thm:nonflat} is proved.

Let us continue with the proof for the second part. With the $C^{1,\sigma}$-regularity of the surface $\partial\{u>0\}\cap B_{\rho r_0}(z)$ at hand, we can infer from the boundary regularity theory for $p$-harmonic functions that $u \in C_{loc}^{1,\eta}(B_{\rho r_0}(z)\cap\{u \geq 0\})$ for some $\eta\in(0,1)$, depending only on $n$, $p$ and $\sigma$. In particular, arguing similarly as in the derivation of \eqref{eq:msr-u0}, we obtain that  
$$
\liminf_{r\to 0} \frac{\Delta_p u^+(B_r(z))}{r^{n-1}} > 0. 
$$ 
Thus, we have proved that 
$$
\partial\{u>0\}\cap \Omega\setminus E \subset \bigcup_{k,\ell=1}^\infty F_{k,\ell},
$$
where $F_{k,\ell}$ consists of all points $z\in \partial\{u>0\}\cap\Omega$ with $\dist(z,\partial\Omega)\geq\frac{1}{\ell}$ such that $\Delta_p u^+(B_r(z)) \geq \frac{1}{k} r^{n-1}$ for any $r\in(0,\frac{1}{2\ell})$. To conclude the proof, it suffices to show that 
$$
H^{n-1} ( F_{k,\ell} ) < \infty,
$$ 
for any pair $(k,\ell)$ of positive integers. However, this claim can be verified by a purely measure-theoretic argument, which has already been shown in the proof of Lemma \ref{lem:msr-deg}. For this reason, we shall not repeat the details, and finish the proof here.



\appendix


\section{Supplementary Detail for Theorem \ref{thm:flat}}\label{appendix:flat}

In what follows, we shall supplement some detail for the proof of Theorem \ref{thm:flat}. We contain several lemmas from \cite{LN1,LN2} with minor modifications made in order to fit the notion of viscosity solutions in Definition \ref{defn:visc}. 

Let us begin with some definition. Let $\cS^n\subset \R^{n^2}$ be the space of all real symmetric $n\times n$ matrices, and $\cS_p^n\subset \cS^n$ be the subspace consisting of all $A\in \cS^n$ whose eigenvalues are bounded above by $\max\{p-1,1\}$ and below by $\min\{p-1,1\}$. Given $p,q\in(1,\infty)$, define $\cP_{p,q}:\cS^n\to\R$ by
$$
\cP_{p,q}(M) = \inf_{A\in \cS_p\cup \cS_q} \tr(AM).
$$

\begin{lem}[Essentially due to {\cite[Lemma 3.5]{LN1}}]\label{lem:flat}
Let $D\subset\R^n$ be a bounded domain, and $\phi\in C^2(D)$ be a positive function such that $\| D \phi \|_{L^\infty(D)} \leq \frac{1}{2}$. Suppose that 
$$
\phi \cP_{p,q} (D^2 \phi) \geq 50\max\{p,q\}n |D \phi|^2\quad\text{in }D,
$$
Let $\Omega\subset\R^n$ be an open set such that $\overline{\bigcup_{x\in D} B_{\phi(x)}(x)} \subset\Omega$, and define $v: D\to \Omega$ by 
$$
v(x) = \sup_{B_{\phi(x)}(x)} u.
$$
Then $v\in C(D)$. Moreover, if $u^+\in W_{loc}^{1,p}(\Omega)$ is $p$-harmonic in $\{u>0\}$ and $u^- \in W_{loc}^{1,q}(\Omega)$ is $q$-harmonic in $\mbox{int}(\{u\leq 0\})$, then $v^+\in W_{loc}^{1,p}(D)$ is a weak $p$-subsolution in $\{ v> 0\}$ and $v^- \in W_{loc}^{1,q}(D)$ is a weak $q$-subsolution in $\mbox{int}(\{ u\leq 0\} )$. 
\end{lem}

\begin{proof}
The assertion that $v\in C(D)$, $v^+ \in W_{loc}^{1,p}(D)$ and $v^-\in W_{loc}^{1,q}(D)$ follows immediately from the assumption $\phi\in C^2(D)$ and $u\in C(\Omega)$. Denote by $\cP_p$ the functional $\cP_{p,p}$. Then $\min\{\cP_p(M),\cP_q(M)\} \geq \cP_{p,q}(M)$ for all $M\in\cS^n$. Hence, it follows from the assumption on $\phi$ that 
$$
\phi \cP_{p_\pm} (D^2 \phi) \geq 50 p_\pm n |D \phi|^2\quad\text{in }D,
$$
where $p_+ = p$ and $p_- = q$. Therefore, one may follow the line of the proof of \cite[Lemma 3.5]{LN1} to derive our conclusion. Note that the argument in the cited reference actually proves that if $u^+$ is $p$-harmonic in $\{ u >0\}$, then $v^+$ is weak $p$-subsolution in $\{ v >0\}$, regardless of the behavior of $u^-$, once the last differential inequality is satisfied. We omit the detail. 
\end{proof}

We shall need some family $\{\phi_t\}_{0\leq t\leq 1}$ of variables radii to construct a continuous family of ``subsolutions''. The following family will be used to prove that Lipschitz free boundaries are $C^{1,\sigma}$.

\begin{lem}[Essentially due to {\cite[Lemma 3.22]{LN1}}]\label{lem:fam1}
Let $\rho \in (0,\frac{1}{100})$, $\gamma\in(0,\frac{1}{2})$ be given. There exist a constant $h\in(0,1]$, depending only on $n$, $p$, $q$ and $\rho$, and a family $\{\phi_t\}_{0\leq t\leq }\subset C^2(B_2\setminus B_\rho(\frac{1}{8}e_n))$ such that 
$$
\begin{cases} 
 \phi_t = 1 &\text{on }\overline{B_2}\setminus B_{1/2} \\
\phi_t \geq 1+h\gamma t &\text{on }B_{1/16}, \\ 
1\leq \phi_t \leq 1+ t\gamma, \quad |D \phi_t| \leq \gamma t & \text{on } \overline{B_2} \setminus \overline{B_\rho(8^{-1}e_n)},\\ 
\phi_t \cP_{p,q} (D^2 \phi) \geq 50\max\{p,q\}n |D \phi_t|^2 &\text{on }\overline{B_2}\setminus B_\rho (8^{-1}e_n).
\end{cases}
$$
\end{lem} 

\begin{proof}
The proof is the same as \cite[Lemma 3.22]{LN1}. The only difference is in the last differential inequality, since we (may) have $p\neq q$. However, this does not contribute to any major modification. More exactly, we can consider the function $f:\R^n\setminus \{0\} \to (0,\infty)$ defined by $f(x) = |x|^{-2N}$, with a sufficiently large $N$ such that $N \geq c\max\{p,q,\frac{1}{p-1},\frac{1}{q-1}\}$, where $c>0$ is a constant depending only on $n$. With such an $f$, one may easily check that $\phi_t = 1+ \frac{t\gamma}{50\max\{p,q\}n} \hat f$, for $0\leq t\leq 1$, satisfies the assertion of this lemma, where $\hat f$ is defined as in the proof of \cite[Lemma 3.22]{LN1}. We skip the detail. 
\end{proof} 

We need another family of variable radii to prove that flat free boundaries are Lipschitz. In what follows, we shall write $x = (x',x_n) \in \R^n$, $Q_{r,s}(x) = \{ (y',y_n)\in\R^n: |y' - x'| < r, |y_n - x_n| <s\}$ and $Q_{r,s} = Q_{r,s}(0)$. 

\begin{lem}[Essentially due to {\cite[Lemma 3.7]{LN2}}]\label{lem:fam2}
Let $\Lambda = \{(x',x_n) \in \R^n: x_n = \lambda (x')\}$ where $\lambda:\R^{n-1}\to \R$ is a function satisfying $\lambda (0) = 0$ and $\| \lambda \|_{\mbox{Lip}(\R^{n-1})} \leq M$ for some constant $M\geq 1$. Given $h\in (0,\frac{1}{10^3})$, let $\Lambda (h) = \{(x',x_n) : |x_n - \lambda (x') < h\} \cap \overline{Q_{2,8M}}$. If $\beta \in (0,1)$, then there exist a family $\{\phi_t\}_{0\leq t\leq 1} \subset C^2(\Lambda(h))$ and constants $c\geq 1$, $h_0 > 0$, both depending only on $n$, $p$, $q$, $M$ and $\beta$, such that the following holds: there exists a constant $\mu\in(0,2]$, depending only on $n$, $p$ and $q$, such that for $x\in h\in (0,h_0]$, 
$$
\begin{cases}
\phi_t = 1 & \text{on } \Lambda (h) \setminus Q_{1-2h^{1-\beta},4M},\\
\phi_t \geq 1 + \mu t - ct h^\beta & \text{on }\Lambda (h) \cap Q_{1-100 h^{1-\beta}, 4M}, \\  
1 \leq\phi_t \leq 1 + \mu t , \quad |D \phi_t|\leq ct h^{\beta -1 }& \text{on }\Lambda (h), \\
\phi_t \cP_{p,q}( D^2 \phi) \geq 50 \max\{p,q\} n |D \phi_t|^2 & \text{on } \Lambda(h), \\
\end{cases}
$$
\end{lem}

\begin{proof}
The proof repeats that of \cite[Lemma 3.7]{LN2}, along with some minor modifications as in the proof of Lemma \ref{lem:fam1}. Thus, we skip the detail. 
\end{proof} 

With the family $\{\phi_t\}_{0\leq t\leq 1}$ of variable radii, we can construct a continuous family $\{v_t\}_{0\leq t\leq 1}$ of deformation as 
$$
v_t (x) = \sup_{B_{\phi_t(x)}(x)} u. 
$$
We can use this family to transfer some strong property of $u$, such as full monotonicity, or $\e$-monotonicity, etc., from a neighborhood of $\{u>0\}\cup \mbox{int}(\{u\leq 0\})$ to its free boundary $\partial\{u>0\}$ in a smaller neighborhood. This is done with the comparison principle and boundary Harnack inequality; the former is available from the definition of viscosity solutions, and the latter is a property of positive weak $p$- (or $q$-) harmonic functions for any $1<p<\infty$. Since $\phi_t$ is constructed such that $\{v_t>0\}\subset \{u>0\}$, the comparison is only performed in one phase (i.e., either $\{u>0\}$ or $\mbox{int}(\{u\leq 0\})$), which is why it does not matter even if $u$ or $v_t$ solves a different operator in the other phase. 

\section*{Declarations}

\noindent {\bf  Data availability statement:} All data needed are contained in the manuscript.

\medskip
\noindent {\bf  Funding and/or Conflicts of interests/Competing interests:} The authors declare that there are no financial, competing or conflict of interests.



\begin{thebibliography}{999999}

\bibitem[AF94]{Acerbi-Fusco}
E. Acerbi and N. Fusco,
{\it A transmission problem in the calculus of variations},
Calc. Var. {\bf 2} (1994), 1--16. 

\bibitem[ACF84]{ACF}
H. W. Alt, L. A. Caffarelli and A. Friedman,
{\it Variational problems with two phases and their free boundaries},
Trans. Amer. Math. Soc. {\bf 282} (1984), 431--461. 

\bibitem[AM12]{AM}
J. Andersson and H. Mikayelyan,
{\it The zero level set for a certain weak solution, with applications to the Bellman equations},
Trans. Amer. Math. Soc. {\bf 365} (2013), 2297--2316. 

\bibitem[BCM18]{BCM18}
P. Baroni, M. Colombo and G. Mingione, 
{\it Regularity for general functionals with double phase},
Calc. Var. {\bf 57} (2018), 1--48. 

\bibitem[Caf87]{C1} 
L. Caffarelli, 
{\it A Harnack inequality approach to the regularity of free boundaries. Part I, Lipschitz free boundaries are $C^{1,\alpha}$},
Rev. Mat. Iberoam. {\bf 3} (1987), 139--162.

\bibitem[Caf89]{C2}
L. Caffarelli, 
{\it A Harnack inequality approach to the regularity of free boundaries. II. Flat free boundaries are Lipschitz},
Comm. Pure Appl. Math. {\bf 42} (1989), 55--78.

\bibitem[Caf88]{C3}
L. Caffarelli, 
{\it A Harnack inequality approach to the regularity of free boundaries. III. Existence theory, compactness, and dependence on $X$}, 
Ann. Sc. Norm. Super. Pisa Cl. Sci. {\bf 15} (1988), 583--602.

\bibitem[CM15a]{Colombo-Mingione-1}
M. Colombo and G. Mingione
{\it Regularity for double phase variational problems},
Arch. for Rational Mech. Anal. {\bf 215} (2015), 443--496.

\bibitem[CM15b]{Colombo-Mingione-2}
M. Colombo and G. Mingione
{\it Bounded minimizers of double phase variational integrals},
Arch. for Rational Mech. Anal. {\bf 218} (2015), 219--273.

\bibitem[CO19]{Chipot}
M. Chipot and H. B. de Oliveira,
{\it Some results on the $p(u)$-Laplacian problem},
Math. Ann. {\bf 375} (2019), 283--306.

\bibitem[DP05]{DP}
D. Danielli and A. Petrosyan,
{\it A minimum problem with free boundary for a degenerate quasilinear operator},
Calc. Var. {\bf 23} (2005), 97--124.

\bibitem[DK18]{DK}
S. Dipierro and A. L. Karakhanyan,
{\it Stratification of free boundary points for a two-phase variational problem},
Adv. Math. {\bf 328} (2018), 40--81.

\bibitem[Fel97]{Fel1}
M. Feldman,
{\it Regularity for nonisotropic two-phase problems with Lipschitz free boundaries},
Differential Integral Equations {\bf 10} (1997), 1171--1179.

\bibitem[Fel01]{Fel2}
M. Feldman,
{\it Regularity of Lipschitz free boundaries in two-phase problems for fully nonlinear elliptic equation},
Indiana Univ. Math. J. {\bf 50} (2001), 1171--1200.

\bibitem[Giu03]{Giu}
E. Giusti,
{\it Direct Methods in the Calculus of Variations},
World Scientific Publishing Co., Inc., River Edge, NJ, 2003.

\bibitem[KLS17]{Kiam-Kim-Sh1}
S. Kim, K.-A. Lee and H. Shahgholian,
{\it An elliptic free boundary arising from the jump of conductivity},
Nonlinear Anal. {\bf 161} (2017), 1--29.

\bibitem[KLS19]{Kiam-Kim-Sh2}
S. Kim, K.-A. Lee and H. Shahgholian,
{\it Nodal Sets for “broken” quasilinear PDEs},
Indiana Univ. Math. J. {\bf 68} (2019), no.4, 1113--1148.

\bibitem[LW17]{LW17}
C. Lederman and N. Wolanski,
{\it Weak solutions and regularity of the interface in an inhomogeneous free boundary problem for the $p(x)$-Laplacian},
Interfaces Free Bound. {\bf 19} (2017), no.2, 201--241.

\bibitem[LW19]{LW19}
C. Lederman and N. Wolanski,
{\it Inhomogeneous minimization problems for the $p(x)$-Laplacian},
J. Math. Anal. Appl. {\bf 475} (2019), 423--463. 

\bibitem[LN10]{LN1}
J. Lewis and K. Nystr\"om, 
{\it Regularity of Lipschitz free boundaries in two-phase problems for the $p$-Laplace operator},
Adv. Math. {\bf 225} (2010), 2565--2597.

\bibitem[LN12]{LN2}
J. Lewis and K. Nystr\"om, 
{\it Regularity of flat free boundaries in two-phase problems for the $p$-Laplace operator},
Ann. I. H. Poincar\'e {\bf 29} (2012), 83--108. 

\bibitem[Lin19]{L}
P. Lindqvist,
{\it Notes on the Stationary $p$-Laplace Equation},
Springer International Publishing, 2019.

\bibitem[Ruz00]{Ruzicka} 
M. Ruzicka,
{\it Electrorheological Fluids: Modeling and Mathematical Theory}, 
Springer-Verlag, Berlin, 2000.


\bibitem[Zhi86]{Zhikov1} 
V. V. Zhikov, {\it Averaging of functionals of the calculus of variations and elasticity theory}, 
Izv. Akad. Nauk SSSR Ser. Mat. {\bf 50} (1986), no. 4, 675--710.

\bibitem[Zhi95]{Zhikov2}
V. V. Zhikov, {\it On Lavrentiev’s phenomenon}, 
Russian J. Math. Phys. {\bf 3} (1995), no. 2, 249--269.



\end{thebibliography}
\end{document}